\documentclass{article}
\usepackage[american]{babel}
\usepackage{amsfonts,amsmath,amssymb}
\input{liemacs.sty}

\newcommand{\out}{\mathop{{\rm out}}\nolimits}
\def\trile{\trianglelefteq}

\def\cyc{{\rm cyc.}}
\def\half{\frac{1}{2}}
\def\k{\fk}
\begin{document}
\title{The second cohomology of  current algebras of general Lie algebras} 

\author{Karl-Hermann Neeb\\
        Fachbereich Mathematik\\
        Technische Universit\"at Darmstadt\\
        Schlossgartenstr. 7\\
        64285 Darmstadt\\
        Germany\\
        neeb@mathematik.tu-darmstadt.de\\
\and
     Friedrich Wagemann \\
     Laboratoire de Math\'ematiques Jean Leray\\
     Facult\'e des Sciences et Techniques\\
     Universit\'e de Nantes\\
     2, rue de la Houssini\`ere\\
     44322 Nantes cedex 3\\
     France\\
%     tel.: ++33.2.51.12.59.57\\
     wagemann@math.univ-nantes.fr}

\maketitle

\begin{abstract}
Let $A$ be a unital commutative associative algebra over a field of 
characteristic zero, $\k$ be a Lie algebra, and 
$\z$ a vector space, considered as a trivial module of the Lie algebra 
$\g := A \otimes \k$. In this paper, we give a 
description of the cohomology space $H^2(\g,\z)$ 
in terms of well accessible data associated to $A$ and $\k$. 
We also discuss the topological situation, where 
$A$ and $\k$ are locally convex algebras. \\


AMS-Class: 17B56, 17B65 
\end{abstract}

\section*{Introduction} 

Let $A$ be a unital commutative associative algebra over a field $\K$ with 
$2 \in \K^\times$ and $\k$ be a $\K$-Lie algebra.
Then the tensor product $\g := A \otimes \k$ is a Lie algebra with 
respect to the bracket 
$$ [a \otimes x, a' \otimes x'] := aa' \otimes [x,x']. $$

Let $\z$ be a vector space, considered as a trivial $\g$-module. The main point of the present paper is to give a 
description of the set $H^2(\g,\z)$ of cohomology classes of $\z$-valued 
$2$-cocycles on the Lie algebra $\g$ in 
terms of data associated to $A$ and $\k$ which is as explicit as possible.

We consider $\z$-valued $2$-cochains on $\g$ as 
linear functions $f \: \Lambda^2(\g) \to \z$. Such a function is a 
{\it $2$-cocycle} if and only if it vanishes on the subspace 
$B_2(\g)$ of $2$-boundaries, which is the image of the linear map 
$$\partial \: \Lambda^3(\g) \to \Lambda^2(\g), \quad 
x \wedge y \wedge z \mapsto [x,y] \wedge z + [y,z] \wedge x + [z,x] \wedge y. $$
In view of the Jacobi identity, $B_2(\g)$ is contained in the subspace 
$Z_2(\g)$ of $2$-cycles, i.e., the kernel of the linear map 
$b_\g \: \Lambda^2(\g) \to \g, x \wedge y \mapsto [x,y].$
The quotient space 
$$ H_2(\g) := Z_2(\g)/B_2(\g) $$
is the {\it second homology space of $\g$}. 

A $2$-cocycle $f$ is a {\it coboundary} if it is of the form 
$f(x,y) = d_\g\ell(x,y) := - \ell([x,y])$ for some linear map $\ell \: \g \to \z$. 
We write $B^2(\g,\z)$ for the set of $2$-coboundaries and 
$Z^2(\g,\z)$ for the set of $2$-cocycles. This 
means that a coboundary vanishes on $Z_2(\g)$. If, conversely, a $2$-cocycle vanishes on 
$Z_2(\g)$, then there exists a linear map $\alpha \: \im(b_\g) = [\g,\g]  \to \z$
with $f = -b_\g^*\alpha$, and any linear extension $\ell$ of $\alpha$ to all of $\g$ 
yields $f = d_\g \ell$. This leads to the following description of the second 
$\z$-valued cohomology group 
$$ H^2(\g,\z) := Z^2(\g,\z)/B^2(\g,\z) \cong \Lin(H_2(\g),\z)\into \Lin(Z_2(\g),\z). $$

{}From this picture, it is clear that we obtain a good description of 
$H^2(\g,\z)$ if we have an accessible description of the space $Z_2(\g)$ 
and its subspace $B_2(\g)$, hence of the quotient space $H_2(\g)$. 
Our goal is a description of this space and the cocycles in terms of 
accessible data attached to the commutative algebra $A$ and the 
Lie algebra $\k$. For a $\K$-vector space $V$, we identify the second exterior power 
$\Lambda^2(V)$, resp., the second symmetric power $S^2(V)$ with the corresponding 
subspaces of $V \otimes V$. Accordingly, we put
$$ x \wedge y := \frac{1}{2}(x \otimes y -y \otimes x) \in \Lambda^2(V)
\quad \mbox{ and } \quad 
 x \vee y := \frac{1}{2}(x \otimes y +y \otimes x) \in S^2(V) $$
and obtain $V \otimes V = \Lambda^2(V) \oplus S^2(V)$. 
For the commutative algebra $A$, we have a natural decomposition 
$S^2(A) \cong A \vee \1 \oplus I_A,$ where 
$I_A \subeq S^2(A)$ is the kernel of the multiplication map $S^2(A) \to A$.
The first step, carried out in Section $2$, is to show that by identifying 
$A$ with $A \vee \1 \subeq S^2(A)$, we obtain a linear isomorphism 
\begin{eqnarray}
  \label{eq:1}
P = (p_1, p_2, p_3) 
\:  \Lambda^2(\g) \to (\Lambda^2(A) \otimes S^2(\k)) \oplus (A \otimes \Lambda^2(\k))
\oplus (I_A \otimes \Lambda^2(\k)), 
\end{eqnarray}
restricting to a linear isomorphism 
\begin{eqnarray}
  \label{eq:2}
Z_2(\g) \to (\Lambda^2(A) \otimes S^2(\k)) \oplus (A \otimes Z_2(\k))
\oplus (I_A \otimes \Lambda^2(\k)).   
\end{eqnarray}
Now each alternating map $f \: \Lambda^2(\g) \to \z$ is represented by three maps 
\begin{eqnarray}
  \label{eq:3}
f_1 \: \Lambda^2(A) \otimes S^2(\k) \to \z, \quad 
f_2 \: A \otimes \Lambda^2(\k) \to \z, 
\quad \mbox{ and } \quad f_3 \: I_A \otimes \Lambda^2(\k) \to \z, 
\end{eqnarray}
determined by $f = \sum_{j=1}^3 f_j \circ p_j$ 
in the sense of (\ref{eq:1}). Since two cocycles define the same cohomology class if and only if they 
coincide on the subspace $Z_2(\g)$ of $\Lambda^2(\g)$, any cohomology class 
$[f] \in H^2(\g,\z)$ is represented by the triple $(f_1, f_2', f_3)$, where 
$f_2' := f_2 \res_{A \otimes Z_2(\k)}$. 
Conversely, three linear maps $f_1, f_2$ and $f_3$ as in (\ref{eq:3}) 
define a cocycle if and only 
if $f := \sum_{j =1}^3 f_j \circ p_j$ vanishes on $B_2(\g)$. 
The main result of the present paper is Theorem~\ref{theo2} which makes 
this condition more explicit as follows: 

(a) The alternating linear map $\tilde f_1 \: A \times A \to \Sym^2(\k,\z)$ defined by 
$\tilde f_1(a,b)(x,y) := f_1(a \wedge b \otimes x \vee y)$ 
has values in the set $\Sym^2(\k,\z)^\k$ of invariant symmetric bilinear maps and 
$f_1$ vanishes on $T_0(A) \otimes (\k \vee \k')$, where 
$$ T_0(A) := \spann \{ ab\wedge c +bc\wedge a + ca\wedge b-  abc \wedge \1 \: 
a,b,c \in A\} $$
and $\k' := [\k,\k]$ denotes the commutator algebra of $\k$. 

(b) For the map $\tilde f_2 \: A \to \Alt^2(\k,\z)$ 
defined by 
$\tilde f_2(a)(x,y) := f_2(a \vee \1 \otimes x \wedge y),$ 
we have 
$$ d_\k(\tilde f_2(a))(x,y,z) = - \tilde f_2(a)(\partial(x \wedge y \wedge z)) 
=  \tilde f_1(a,\1)([x,y],z) 
\quad \mbox{ for all} \quad a \in A, x,y,z \in \k,$$
with the Lie algebra differential $d_\k \: C^2(\k,\z) = \Alt^2(\k,\z) \to Z^3(\k,\z)$. 

(c) $f_3$ vanishes on $I_A \otimes (\k \times \k')$. 

Note that these conditions imply that the two maps 
$f_1 \oplus f_2$ and $f_3$ are also cocycles, 
whereas $f_1$ and $f_2$ are cocycles if and only if $f_1$ vanishes on 
$(A \wedge \1) \otimes (\k \vee \k')$, which, in view of (b), means that 
$\tilde f_2(A)$ vanishes on $B_2(\k)$, i.e., $\tilde f_2$ has values in the 
space $Z^2(\k,\z)$ of $\z$-valued $2$-cocycles on $\k$. 
Cocycles of the form $f_1 \oplus f_2$, where $f_1$ and $f_2$ are not cocycles, are called 
{\it coupled}. All coboundaries are of the form $f = f_2$ ($f_1 = f_3 = 0$), 
so that the cohomology class of a coupled 
cocycle contains only coupled cocycles.  

We show that $\g$ possesses  non-zero coupled $2$-cocycles
 if and only if the image of the universal derivation $d_A \: A \to \Omega^1(A)$ is non-trivial 
and $\fk$ possesses a symmetric invariant bilinear form $\kappa$ for which the 
$3$-cocycle $\Gamma(\kappa)(x,y,z):=\kappa([x,y],z)$ is a non-zero coboundary. 
The map $\Gamma \: \Sym^2(\k)^\k \to Z^3(\k)$ is called the {\it Koszul map} 
(cf.\ \cite{Kos50}, \S 11; see also \cite{ChE48}, p.113).  
Calling an invariant symmetric bilinear form $\kappa \in \Sym^2(\k)^\k$ {\it exact} 
if $\Gamma(\kappa)$ is a coboundary, this means that $\k$ possesses exact 
invariant bilinear forms $\kappa$ with $\Gamma(\kappa)$ non-zero. 
Note that this is 
not the case if $\k$ is finite-dimensional semisimple, so that there are 
no coupled cocycles in this case. 

Our approach leads us to an exact sequence of the form 
$$ \{0\} \to H^2(\g/\g')_{1,3} \oplus \Lin(A,H^2(\k)) \ssmapright{\Phi} H^2(\g) 
\ssmapright{\Psi} \Lin((\Omega^1(A), d_A(A)), (Z^3(\k)_\Gamma, B^3(\k)_\Gamma)) \to \{0\},$$ 
which is the main result of Section~4. 
Here $H^2(\g/\g')_{1,3}$ denotes the set of alternating bilinear forms 
on $\g/\g' \cong A \otimes \k/\k'$ of the form $f_1 + f_3$, 
and for two pairs $(X,X')$ and $(Y,Y')$ of linear spaces with $X' \subeq X$ and 
$Y' \subeq Y$ we write 
$$ \Lin((X,X'), (Y,Y')) := \{ f \in \Lin(X,Y) \: f(X') \subeq Y'\}, $$
so that we have $\Lin(X,Y) =  \Lin((X,\0), (Y,\0))$, 
and we put $Z^3(\k)_\Gamma := \im(\Gamma) \subeq Z^3(\k)$ and 
$B^3(\k)_\Gamma := B^3(\k) \cap \im(\Gamma).$
{}From the exact sequence, it follows that a crucial part of the description 
of $H^2(\g)$ lies in an understanding of the spaces 
$Z^3(\k)_\Gamma$ and $B^3(\k)_\Gamma$. In an appendix, we show that 
the map $\gamma \: \Sym^2(\k)^\k \to H^3(\k)$, 
induced by the Koszul map $\Gamma$, is part of an exact sequence 
\begin{eqnarray} \label{exseq0} 
&& \{0\} \to H^2(\k) \to H^1(\k,\k^*) \to \Sym^2(\k)^\k \ssmapright{\gamma} H^3(\k) 
\to H^2(\k,\k^*) \to H^1(\k,\Sym^2(\k)),  
\end{eqnarray}
which implies that for the space $\Sym^2(\k)^\k_{\rm ex} = \ker \gamma$ 
of exact invariant forms we have 
$$ \Sym^2(\k)^\k_{\rm ex} \cong H^1(\k,\k^*)/H^2(\k)
\quad \mbox{ and } \quad 
\im(\gamma) \cong \Sym^2(\k)^\k/\Sym^2(\k)^\k_{\rm ex}
\cong Z^3(\k)_\Gamma/B^3(\k)_\Gamma.$$ 

In Section~\ref{sec5}, we give an example of a non-trivial 
coupled $2$-cocycle and in Section~\ref{sec6} we explain 
how our results can be used for the analysis of continuous cocycles 
if $\K \in \{\R,\C\}$ and $A$ and $\k$ are locally convex spaces with continuous algebra 
structures. Then $\g = A \otimes \k$ carries the structure of a locally 
convex Lie algebra, and we are interested in the space 
$H_c^2(\g,\z)$ of cohomology classes of continuous $2$-cocycles with values in a 
locally convex space $\z$ modulo those coboundaries coming from continuous linear maps 
$\g \to \z$. The main difficulty in applying the algebraic results in the topological 
context with an infinite dimensional Lie algebra $\k$ 
 is the possible discontinuity of a linear map $h: \g \to \z$ bounding an algebraically 
trivial $2$-cocycle.  

If $\fk$ is a finite dimensional semi-simple Lie algebra and $A$ a 
topological algebra, then the continuous second cohomology 
space $H_c^2(\g,\K)$ has been determined in \cite{Mai} as $Z^1_c(A,\K) \otimes 
\Sym^2(\k,\K)^\k$, where $Z^1_c(A,\K)$ denotes the space of continuous $\K$-valued 
cyclic $1$-cocycles on $A$ (see \cite{KL82} for the algebraic case). 
As any exact form vanishes on a semi-simple 
Lie algebra, there are no coupled cocycles in this case. 
%An interesting open problem is to determine
% conditions on Lie algebras for the existence of exact forms.  

The main previous contributions to the investigations of $H_2(\g)$ for $\g = A \otimes \k$ and 
arbitrary $\k$ and $A$ are the articles by Haddi \cite{Haddi} and Zusmanovich \cite{Sus}.
 Both offer a description of $H_2(\g)$
 in terms of (sub- or quotient) spaces. Haddi \cite{Haddi} uses the projection $s_2:H_2(\fg)\to H_2((\Lambda^*(\fg)_{\fk},\partial))$ of $H_2(\fg)$ to the homology of the quotient complex of $\fk$-coinvariants and computes kernel and cokernel of this map. The cokernel of $s_2$ is isomorphic to 
$\g'/[\g,\g']$, and the kernel is isomorphic to 
$(A\otimes \oline H_2(\fk,\fk'))\oplus D(A,\fk,\fk')$, 
$\oline H_2(\fk,\fk')$ is the kernel of 
the projection $H_2(\fk)\to H_2(\fk/\fk')$ (the subspace of {\it essential 
homology}), 
and $D(A,\fk,\fk')$ is the subspace of $H_2(\fg)$ 
generated by cycles of the form $ax\wedge y+ay\wedge x$ for $x$ or $y\in\fk'$ and $a\in A$, 
which lies in $\Lambda^2(A) \otimes S^2(\k)$ (in our notation). 
Furthermore he uses a non-canonical splitting to identify the homology of the 
coinvariants $H_2((\Lambda^*(\fg)_{\fk},\partial))$ 
with $(\Omega^1(A)/d_A(A)\otimes B_{\fk,\fk'})
\oplus\Lambda^2(\g/\g')$, where $B_{\fk,\fk'}$ 
is the image of $\fk\vee\fk'$ in the space of $\fk$-coinvariants of symmetric $2$-tensors on $\fk$. He thus obtains an exact sequence
$$0\to (A\otimes \oline H_2(\fk,\fk'))\oplus D(A,\fk,\fk')\to H_2(\fg)
\ssmapright{s_2} 
\big(\Omega^1(A)/d_A(A)\otimes B_{\fk,\fk'}\big)\oplus\Lambda^2(\g/\g')\to \g'/[\g,\g'] 
\to 0.$$ 
It is instructive to compare this sequence with our exact cohomology sequence 
described above. 

Zusmanovich \cite{Sus} uses as extra data a free presentation of $\fk$ and deduces one of $\g$. He describes the subspace of essential homology 
$\oline H_2(\fg,\fg')$ by the Hopf formula in terms of the presentation. 
In this way, he identifies the different terms in the exact sequence given by the $5$-term
 exact sequence of the Hochschild--Serre spectral sequence for the subalgebra $\fg'\subset\fg$ 
(using non-canonical splittings). His description yields
$$H_2(\fg)\simeq (A \otimes H_2(\fk))\,\oplus\,(\Omega^1(A)/d_A(A) \otimes B(\fk))\oplus\,
(\Lambda^2(\fk/\fk')\otimes I_A) \oplus\,(S^2(\fk/\fk')\otimes\,T(A)),$$
where $B(\fk)$ is the space of $\fk$-coinvariants in $S^2(\fk)$, 
and $T(A)\subset\Lambda^2(A)$ is 
spanned by the elements $ab\wedge c+ca\wedge b+bc\wedge a$ for $a,b,c\in A$. 
 
The main advantage of our approach is that is does not require any auxiliary data 
and provides a quite explicit description of cocycles representing the different 
types of 
cohomology classes. In particular, this direct approach leads us to the 
interesting new class of coupled cocycles. In subsequent work, we plan to 
use the methods developed in \cite{Ne02} to study 
global central extensions of Lie groups $G$ whose Lie algebras are of the form 
$\g = A \otimes \k$ defined by coupled Lie algebra cocycles. 
For algebras of the type $A = C^\infty_c(M,\R)$, i.e., compactly supported smooth 
functions on a manifold $M$, this has been carried out in \cite{MN03} and 
\cite{Ne04}. 

{\bf Thanks:} We are grateful to M.~Bordemann for a stimulating email exchange and for 
pointing out the relation to the exact sequence (\ref{exseq0}), 
part of which is due to him. We also thank the referee for an extremely valuable 
report and in particular for pointing out several references and some inaccuracies in 
a previous version. 

\subsection*{Notation} 

In the following, we write elements of $\g = A \otimes \fk$ simply as 
$ax := a \otimes x$ to simplify notation. Elements of $A$ are mostly 
denoted $a,b,c,\ldots$ or $a,a',a'',\ldots$ and elements 
of $\fk$ are denoted $x,y,z,\ldots$ or $x,x',x'',\ldots$. 
We write $\k' := [\k,\k]$ for the commutator algebra of $\k$ and observe that 
$\g' = A \otimes \k'$ is the commutator algebra of $\g$. 

We also write 
$C^p(\g) := C^p(\g,\K)$, 
$Z^p(\g) := Z^p(\g,\K)$, 
$B^p(\g) := B^p(\g,\K)$, and  
$H^p(\g) := H^p(\g,\K)$ for the spaces of Lie algebra $p$-cochains, cocycles, coboundaries 
and cohomology classes with values in the trivial module $\K$. 
We write $\Sym^2(\k,\z)$ for the space of $\z$-valued 
symmetric bilinear maps on $\k$ and put $\Sym^2(\k) := \Sym^2(\k,\K)$. 
Accordingly, we write $\Alt^2(\k,\z)$ for the set of $\z$-valued alternating bilinear maps. 

\section{Several approaches to the universal differential module of $A$} \label{sec1}

In this section, we review different constructions of the universal differential 
module $\Omega^1(A)$. The relationship between these constructions will play a crucial 
role in the following. 

An important object attached to the algebra $A$ is 
its universal differential module $\Omega^1(A)$. This is an $A$-module 
with a derivation $d_A \: A \to \Omega^1(A)$ which is
universal in the sense that for any other $A$-module $M$ 
and any derivation $D \: A \to M$, there exists a unique 
module morphism $\alpha \: \Omega^1(A) \to M$ of $A$-modules 
with $D= \alpha\circ d_A$. From its universal property, it is easy to derive that 
the universal differential module is unique up to isomorphism, 
but there are many realizations, looking at first sight quite differently. 

Let $\mu_A \: A \otimes A \to A, a \otimes b \mapsto ab$ denote the multiplication 
of $A$. Then $\mu_A$ is an algebra morphism, so that $J_A := \ker \mu_A$ is an ideal 
of the commutative algebra $A \otimes A$. 
{}From the $A$-module structure on $A \otimes A$, given by 
$a.(b \otimes c) := ab \otimes c$, we thus derive 
an $A$-module structure on the quotient space $J_A/J_A^2$, 
which also is a (non-unital) commutative algebra. 
Let $[x]$ denote the image of $x \in J_A$ in $J_A/J_A^2$. Then 
$$ D \: A \to J_A/J_A^2, \quad a \mapsto [\1 \otimes a- a \otimes \1] $$
is a derivation and it is not hard to verify that 
$(J_A/J_A^2, D)$ has the universal property of $(\Omega^1(A),d_A)$ 
(cf.\ \cite{Bou90}, Ch.~III, \S 10.11). We obviously have the direct decomposition 
$A \otimes A =(A \otimes \1) \oplus J_A$, where the projection onto the 
subspace $J_A$ is given by 
$$p \: A \otimes A \to J_A, \quad 
a \otimes b \mapsto a \otimes b - ab \otimes \1.$$ 
This implies that 
$$ J_A 
= \spann \{ a \otimes b - ab \otimes \1 \: a,b \in A \} 
= (A \otimes \1).\spann \{ \1 \otimes b - b \otimes \1 \: b \in A \}, $$
and thus 
\begin{eqnarray}   \label{eq:5}
J_A^2 
&=& \spann \{ (a \otimes \1) (\1 \otimes b - b \otimes \1)
(\1 \otimes c - c \otimes \1) \: a,b,c \in A \} \\ 
&=& \spann \{ a \otimes bc - ab \otimes c - ac \otimes b + abc \otimes \1 \: 
a,b,c \in A\}. \notag 
\end{eqnarray}

Another way to construct $\Omega^1(A)$ is by observing that 
each linear map $D \: A \to M$ leads to a linear map 
$\tilde D \: A \otimes A \to M, a \otimes b \mapsto a Db,$ 
and that $D$ is a 
derivation if and only if 
$$ \ker \tilde D \supeq \{ \1 \otimes ab - a \otimes b - b \otimes a \: a,b \in A\}, $$
which implies that $\ker \tilde D$ contains the $A$-submodule 
$$ B_1(A) := \spann \{ c \otimes ab - ca \otimes b - cb \otimes a \: a,b,c \in A\} 
= \spann \{ ab \otimes c + ac \otimes b - a \otimes bc \: a,b,c \in A\},  $$
of $A \otimes A$. The quotient 
$$ HH_1(A) := (A \otimes A)/B_1(A) $$
is called the {\it first Hochschild homology space of $A$}. 
{}From the preceding discussion, it follows that the map 
\begin{eqnarray} \label{hhomega}
HH_1(A) \to \Omega^1(A), \quad [a \otimes b] \mapsto a d_A(b) 
\end{eqnarray} 
is an isomorphism of $A$-modules because the map 
$D \: A \to HH_1(A), a \mapsto [\1 \otimes a]$
is a derivation with the universal property (cf.\ \cite{Lod}, Prop.\ 1.1.10). 
The link between the description of $\Omega^1(A)$ as $HH_1(A)$ and 
$J_A/J_A^2$ is given by the commutative diagram 
$$ 
\begin{matrix} 
A \otimes A & \ssmapright{p} & J_A \\ 
\mapdown{} & & \mapdown{} \\ 
HH_1(A) & \ssmapright{\phi}& J_A/J_A^2 \\
\end{matrix} 
$$
with the isomorphism $\phi([a \otimes b]) = a D(b) = [a \otimes b - ab \otimes \1]$.  
Note that the commutativity of the diagram implies that 
\begin{eqnarray} \label{jb}
J_A^2 = p(B_1(A)). 
\end{eqnarray} 

Let 
$$T(A) := \spann \{ ab\wedge c +bc\wedge a + ca\wedge b \in \Lambda^2(A) \: a,b,c \in A \} $$
denote the image of the subspace $B_1(A) \subeq A \otimes A$ under the quotient map 
$A \otimes A \to \Lambda^2(A), a \otimes b \mapsto a \wedge b$. 
In view of $a d_A(b) + b d_A(a) = d_A(ab)$, the image of the subspace of 
symmetric tensors, which we identify with $S^2(A)$, in 
$\Omega^1(A)$ coincides with $d_A(A)$, so that  (\ref{hhomega}) 
immediately shows that the map 
$$ \Lambda^2(A)/T(A) \cong (A \otimes A)/(S^2(A) + B_1(A)) 
\to \Omega^1(A)/d_A(A), \quad 
[a \wedge b] \mapsto [a d_A(b)] $$
induces a linear isomorphism. 
It is well known that the {\it first cyclic homology space} 
$$ HC_1(A) := \Omega^1(A)/d_A(A) \cong HH_1(A)/[\1 \otimes A] = \Lambda^2(A)/T(A) $$
is of central importance for Lie algebra $2$-cocycles on Lie algebras of the form 
$A \otimes \k$ (cf.\ \cite{KL82}). 

Alternating bilinear maps $f \: A \times A \to \z$ for which the corresponding map 
$\Lambda^2(A) \to \z$ vanishes on $T(A)$ are called {\it cyclic $1$-cocycles}, 
which means that 
$$ f(a,bc) + f(b,ca) + f(c,ab) = 0 \quad \mbox{ for } \quad a,b,c \in A. $$
{}From the above, it follows that the space $Z^1(A,\z)$ of $\z$-valued cyclic $1$-cocycles 
can be identified with 
$$\Lin(HC_1(A), \z) \cong  \{ L \in \Lin(\Omega^1(A),\z) \: d_A(A) \subeq \ker L\}. $$

We define two trilinear maps 
$$ T \: A^3  \to \Lambda^2(A), \quad 
(a,b,c) \mapsto \sum_\cyc ab \wedge c := ab\wedge c +bc\wedge a + ca\wedge b $$
and 
$$ T_0 \: A^3  \to \Lambda^2(A), \quad 
T_0(a,b,c) := T(a,b,c) - abc \wedge \1. $$
We also put $T_0(A) := \spann(\im(T_0))$. 

\begin{Lemma} \label{lem1.1} The map 
$$ \gamma_A \: \Lambda^2(A) \to \Omega^1(A),\quad a \wedge b \mapsto ad_A(b) - b d_A(a) $$
 is surjective and $\ker \gamma_A = T_0(A)$. 
\end{Lemma} 

\begin{Proof} That $\gamma_A$ is surjective follows from 
$$ \gamma_A(a \wedge b + \1 \wedge ab) 
= ad_A(b) - b d_A(a) + d_A(ab) = 2 a d_A(b). $$

For the determination of the kernel of $\gamma_A$, we use 
the realization of $\Omega^1(A)$ as $J_A/J_A^2$. 
In this case, $d_A(a) = [\1 \otimes a - a \otimes \1]$, so that 
$$ \gamma_A(a \wedge b) 
= [a \otimes b - ab \otimes \1  
- b \otimes a + ba \otimes \1] 
= [a \otimes b - b \otimes a]. $$
Therefore the kernel of $\gamma_A$ is the intersection 
$\Lambda^2(A) \cap J_A^2,$
where we consider $\Lambda^2(A)$ as the subspace of skew-symmetric tensors in 
$A \otimes A$. 

Writing $A \otimes A$ as $\Lambda^2(A) \oplus S^2(A)$, the commutativity 
of the multiplication of $A$ shows that
\begin{eqnarray} \label{ja-dec}
J_A = \Lambda^2(A) \oplus I_A \quad \mbox{holds for } \quad 
I_A := J_A \cap S^2(A). 
\end{eqnarray}
Since the flip involution is an algebra isomorphism of $A \otimes A$, we 
have 
$$ \Lambda^2(A) \Lambda^2(A) + S^2(A) S^2(A) \subeq S^2(A)
\quad \mbox{ and } \quad 
\Lambda^2(A) S^2(A) \subeq \Lambda^2(A). $$
This implies that 
$$ \ker \gamma_A = \Lambda^2(A) \cap J_A^2 = I_A \cdot \Lambda^2(A), $$
and that this subspace coincides with the image of $J_A^2$ under the projection 
$$ \alpha \: A \otimes A \to \Lambda^2(A), \quad a \otimes b \mapsto 
a \wedge b = \frac{1}{2}(a \otimes b - b \otimes a). $$
Finally, this leads with (\ref{eq:5}) to 
$$ \ker \gamma_A = \alpha(J_A^2) 
= \spann \{ a \wedge bc - ab \wedge c - ac \wedge b + abc \wedge \1 \: 
a,b,c \in A\} = T_0(A). $$
\end{Proof}

\section{A decomposition of $\Lambda^2(\g)$} \label{sec2}

In this section, we turn to the identification of the space $B_2(\g)$ 
of $2$-coboundaries in $Z_2(\g)$ in terms of our threefold direct sum 
decomposition (\ref{eq:2}). 

{}From 
the universal property of $\Lambda^2(\g)$, we immediately obtain linear maps 
$$p_+ \: \Lambda^2(\g) \to \Lambda^2(A) \otimes S^2(\k), \quad ax \wedge by \mapsto 
a \wedge b \otimes x \vee y $$
and 
$$p_- \: \Lambda^2(\g) \to S^2(A) \otimes \Lambda^2(\k), \quad ax \wedge by 
\mapsto a \vee b \otimes x \wedge y.  $$
We likewise have linear maps 
$$ \sigma_+ \: \Lambda^2(A) \otimes S^2(\k) \to \Lambda^2(\g), 
\quad a \wedge b \otimes x \vee y \to \half(ax \wedge by + ay \wedge bx) $$
and 
$$ \sigma_- \: S^2(A) \otimes \Lambda^2(\k) \to \Lambda^2(\g), 
\quad a \vee b \otimes x \wedge y \to \half(ax \wedge by - ay \wedge bx) $$
satisfying 
$$ p_+ \circ \sigma_+ = \id, \quad p_- \circ \sigma_- = \id \quad \mbox{ and } \quad 
\sigma_+ p_+ + \sigma_- p_- = \id_{\Lambda^2(\g)}. $$
In this sense, we have 
$$ \Lambda^2(\g) \cong \big(\Lambda^2(A) \otimes S^2(\k)\big) \oplus 
\big(S^2(A) \otimes \Lambda^2(\k)\big), $$
and the projections on the two summands are given by $p_\pm$. 

Recall the kernel $J_A$ of the multiplication map $\mu_A \: A \otimes A \to A$. 
The map 
$$ \sigma_A \: A \to S^2(A), \quad a \mapsto a \vee \1 $$
is a section of the multiplication map $\mu_A$, so that we obtain 
a direct sum decomposition 
$$ S^2(A) = (A \vee \1) \oplus I_A \cong A \oplus I_A $$
(cf.\ (\ref{ja-dec})).
In view of this decomposition, we obtain a linear isomorphism 
\begin{eqnarray}
  \label{eq:new}
P = (p_1, p_2, p_3) 
\:  \Lambda^2(\g) \to (\Lambda^2(A) \otimes S^2(\k)) \oplus (A \otimes \Lambda^2(\k))
\oplus (I_A \otimes \Lambda^2(\k)), 
\end{eqnarray}
where the projections $p_1, p_2, p_3$ on the three summands are given by 
$$ p_1(ax \wedge by) = p_+(ax \wedge by) = a \wedge b \otimes x \vee y 
= \frac{1}{2}(ax \wedge by + ay \wedge bx), $$
$$ p_2(ax \wedge by) = ab \otimes x \wedge y, \quad \mbox{ and } \quad 
p_3(ax \wedge by) = (a \vee b - ab \vee \1) \otimes x \wedge y. $$

The following lemma provides the decomposition of $Z_2(\g)$ which is a 
central tool in the following. 
\begin{Lemma} \label{lem2.1} The space $Z_2(\g)$ is adapted to the 
direct sum decomposition of $\Lambda^2(\g)$: 
$$ P(Z_2(\g)) = (\Lambda^2(A) \otimes S^2(\k)) \oplus 
(A \otimes Z_2(\k)) \oplus (I_A \otimes \Lambda^2(\k)). $$
\end{Lemma} 

\begin{Proof} Since $b_\g(ax \wedge by) = ab [x,y]$ is symmetric in $a,b$ and 
alternating in $x,y$, its kernel contains $\Lambda^2(A) \otimes S^2(\k)$. 
The formula for $b_\g$ also shows immediately that 
$I_A \otimes \Lambda^2(\k) \subeq \ker b_\g$, so that 
it remains to observe that 
$$ P(Z_2(\g)) \cap (A \otimes \Lambda^2(\k)) = \ker b_\g \cap (A \otimes \Lambda^2(\k))
= A \otimes Z_2(\k) $$
because 
$b_\g(a \vee \1 \otimes x \wedge y) 
=  \half b_\g(ax \wedge y + x \wedge ay) = a [x,y] = a b_\k(x \wedge y).$
\end{Proof}

In the following, we write $\equiv \mod B_2(\g)$ for congruence of elements 
of $\Lambda^2(\g)$ modulo $B_2(\g)$.  
\begin{Lemma} \label{lem2.2} 
For $a,b,c \in A$ and $x,y,z \in \k$, we have 
$$ p_1(\partial(ax\wedge by \wedge cz)) \equiv p_1(abc x\wedge [y,z])
\equiv -  p_2(\partial(ax\wedge by \wedge cz)) \mod B_2(\g). $$
In particular, $(p_1 + p_2)(B_2(\g)) \subeq B_2(\g).$
\end{Lemma} 

\begin{Proof} From 
\begin{eqnarray*}
\partial(ax \wedge by \wedge cz) &=& ab[x,y] \wedge cz + bc[y,z] \wedge ax + ac [z,x] \wedge by \\
\partial(cx \wedge ay \wedge bz) &=& ac[x,y] \wedge bz + ab[y,z] \wedge cx + bc [z,x] \wedge ay \\
\partial(x \wedge ac y \wedge bz) &=& ac[x,y] \wedge bz + abc[y,z] \wedge x + b [z,x] \wedge acy \\
\partial(cx \wedge y \wedge abz) &=& c[x,y] \wedge abz + ab[y,z] \wedge cx + abc [z,x] \wedge y \\
\partial(bcx \wedge ay \wedge z) &=& abc[x,y] \wedge z + a[y,z] \wedge bcx + bc [z,x] \wedge ay \\
\partial(abcx \wedge y \wedge z) &=& abc[x,y] \wedge z + [y,z] \wedge abcx +abc [z,x] \wedge y, \\
\end{eqnarray*}
we derive 
\begin{eqnarray*}
&&  \partial(ax \wedge by \wedge cz)
+\partial(cx \wedge ay \wedge bz) 
-\partial(x \wedge ac y \wedge bz)
-\partial(cx \wedge y \wedge abz) \\
&& - \partial(bcx \wedge ay \wedge z) 
+ \partial(abcx \wedge y \wedge z) \\
&=& ab[x,y] \wedge cz + bc[y,z] \wedge ax + ac [z,x] \wedge by  
- b[z,x] \wedge ac y - c[x,y] \wedge abz - a [y,z] \wedge bc x \\
&& - abc[y,z] \wedge x + [y,z] \wedge abcx \\
&=& 2 ab \wedge c \otimes [x,y] \vee z 
+ 2 bc \wedge a \otimes [y,z] \vee x + 2 ac \wedge b \otimes [z,x] \vee y - 2 abc \wedge 1 \otimes x \vee [y,z] \\
&=& 2 p_1(\partial(ax \wedge by \wedge cz)) - 2 p_1(abc x \wedge [y,z]). \end{eqnarray*}
This proves the first congruence. 

Note that for $a \in A$ and $x,y,z \in \k$ we have 
$$ \partial(ax \wedge y \wedge z) = a[x,y] \wedge z + [y,z] \wedge ax +a [z,x] \wedge y, $$
which implies that 
\begin{eqnarray} \label{rel1}
a[x,y] \wedge z + a[z,x] \wedge y  \equiv ax \wedge [y,z]  \mod B_2(\g).    
\end{eqnarray}
Summing over all cyclic permutations of $(x,y,z),$ leads to 
\begin{eqnarray}\label{rel2} 
 2\sum_\cyc a[x,y] \wedge z  \equiv \sum_\cyc ax \wedge [y,z] \mod B_2(\g). 
\end{eqnarray}

{}From the relation (\ref{rel1}), we get 
$$ 2 p_1(abc x \wedge [y,z])
= abc[y,z] \wedge x +  abc x \wedge [y,z] 
\equiv \sum_\cyc  abc[y,z] \wedge x = \sum_\cyc  abc[x,y] \wedge z. $$
In view of 
\begin{eqnarray}  
p_2(\partial(ax \wedge by \wedge cz)) 
&=& p_2(ab[x,y]\wedge cz + bc[y,z]\wedge ax + ca[z,x]\wedge by) \notag \\
&=& abc\otimes([x,y]\wedge z + [y,z]\wedge x + [z,x]\wedge y) \notag \\
&=& abc \otimes \partial(x \wedge y \wedge z), \label{rel3}
\end{eqnarray}
relation (\ref{rel2}) yields \pagebreak 
\begin{eqnarray*}
 2 p_2(\partial(ax \wedge by \wedge cz)) 
&=& \sum_\cyc abc[x,y] \wedge z - \sum_\cyc abc z \wedge [x,y]  
\equiv \sum_\cyc abc[x,y] \wedge z - 2 \sum_\cyc abc[x,y] \wedge z \\
&=&  -\sum_\cyc abc[x,y] \wedge z \equiv 
- 2 p_1(\partial(ax \wedge by \wedge cz)). 
\end{eqnarray*}
\end{Proof} 

In view of the preceding lemma, the projection 
$p_1 + p_2$ of $\Lambda^2(\g)$ onto the subspace 
$$ \Lambda^2(A) \otimes S^2(\k) \oplus A \otimes \Lambda^2(\k) $$
preserves $B_2(\g)$. This also implies that 
$\id - p_1 - p_2 = p_3$ preserves $B_2(\g)$, and we derive that 
$$ B_2(\g) = B_2(\g) \cap \big(\Lambda^2(A) \otimes S^2(\k) \oplus A \otimes \Lambda^2(\k)) 
\oplus B_2(\g) \cap (I_A \otimes \Lambda^2(\k)). $$
The following lemma provides refined information. 
\begin{Lemma} \label{lem2.3} 
  \begin{description}
  \item[(1)] $\Lambda^2(A) \otimes \k.S^2(\k) + T_0(A) \otimes \k \vee \k' 
\subeq B_2(\g)$ and $p_1(B_2(\g)) 
= \Lambda^2(A) \otimes \k.S^2(\k) + T(A)\otimes \k\vee\k'$.
  \item[(2)] $p_2(B_2(\g)) = A\otimes B_2(\k)$.
  \item[(3)] $I_A \otimes (\k \wedge \k')= p_3(B_2(\g)) \subeq B_2(\g)$. 
  \end{description} 
\end{Lemma} 

\begin{Proof} (2) follows immediately from formula (\ref{rel3}). 

(1) Recall the identifications 
$x \wedge y = \frac{1}{2}(x \otimes y -y \otimes x)$ 
and $x \vee y = \frac{1}{2}(x \otimes y +y \otimes x)$. 
That $\Lambda^2(A) \otimes \k.S^2(\k)$ is contained in $B_2(\g)$ 
follows immediately from 
\begin{eqnarray*}
&& \partial(ax \wedge by \wedge z) - \partial(bx \wedge ay \wedge z) \\
&=& ab[x,y] \wedge z + b[y,z] \wedge ax + a[z,x] \wedge by 
- ab[x,y] \wedge z - a[y,z] \wedge bx - b[z,x] \wedge ay  \\ 
&=&  b[y,z] \wedge ax + a[z,x] \wedge by  - a[y,z] \wedge bx - b[z,x] \wedge ay  \\ 
&=&  2b \wedge a \otimes [y,z] \vee x + 2 a \wedge b \otimes [z,x] \vee y \\ 
&=& 2 a \wedge b \otimes ([z,x] \vee y - [y,z] \vee x)\\ 
&=& 2 a \wedge b \otimes z.(x \vee y).
\end{eqnarray*}
Therefore the description of $p_1(B_2(\g)) = \im(p_1 \circ \partial)$ follows from 
\begin{eqnarray*}
p_1(\partial(ax \wedge by \wedge cz)) 
&=& p_1(ab[x,y]\wedge cz + bc[y,z]\wedge ax + ca[z,x]\wedge by) \\
&=& ab\wedge c\otimes[x,y]\vee z + bc\wedge a\otimes[y,z]\vee x + ca\wedge b\otimes[z,x]\vee y \\
&\equiv & (ab\wedge c +bc\wedge a + ca\wedge b)\otimes[x,y]\vee z\quad \mod 
\Lambda^2(A)\otimes\k.S^2(\k) \\
&=& T(a,b,c) \otimes[x,y]\vee z. 
\end{eqnarray*}

In (\ref{rel3}), we have seen that 
$$ p_2(\partial(ax \wedge by \wedge cz)) = abc \otimes \partial(x \wedge y \wedge z), $$
and this implies that 
$$ p_2(\partial(abx \wedge y \wedge cz)) = abc \otimes \partial(x \wedge y \wedge z), $$
which leads to 
$$ \partial(ax \wedge by \wedge cz) - \partial(abx \wedge y \wedge cz)
\in \ker p_2. $$
In view of 
$$ T(a,b,c) - T(ab,\1,c)
= T(a,b,c) - (ab \wedge c + c \wedge ab + abc \wedge \1) 
= T(a,b,c) - abc \wedge \1 = T_0(a,b,c) $$
and Lemma~\ref{lem2.2}, the following element is contained in $B_2(\g)$: 
\begin{eqnarray*}
p_1(\partial(ax \wedge by \wedge cz) - \partial(abx \wedge y \wedge cz))
&\in & \big(T(a,b,c) - T(ab,\1,c)\big) \otimes[x,y]\vee z + \Lambda^2(A)\otimes\k.S^2(\k) \\
&\subeq& T_0(a,b,c)\otimes[x,y]\vee z + B_2(\g),
\end{eqnarray*}
and now Lemma~\ref{lem2.2} implies that $T_0(A) \otimes \k' \vee \k \subeq B_2(\g)$.

(3) First we note that 
\begin{eqnarray}
&& p_3(\partial(ax \wedge by \wedge z)) 
=  p_3(ab[x,y] \wedge  z + b[y,z] \wedge ax + a[z,x] \wedge by) \notag \\
&=& (ab \vee \1 - ab \vee \1) \otimes [x,y] \wedge z 
+ (b \vee a - ab \vee \1) \otimes [y,z] \wedge x + (a \vee b - ab \vee \1) \otimes [z,x]\wedge y\notag  \\
&=& (a \vee b - ab \vee \1) \otimes ([y,z] \wedge x + [z,x]\wedge y). \label{eq:sharp}
\end{eqnarray}
Since $p_3$ preserves $B_2(\g)$ (Lemma~\ref{lem2.2}), 
this expression lies in $B_2(\g)$. 
Using the same formula for all cyclic permutations of $x,y,z$ and adding all three terms, 
we see that 
$$ 2(a \vee b - ab \vee \1) \otimes \sum_\cyc [x,y] \wedge z \in B_2(\g). $$ 
This also implies that 
$$ (a \vee b - ab \vee \1) \otimes [x,y] \wedge z 
=  (a \vee b - ab \vee \1) \otimes \Big(\sum_\cyc [x,y] \wedge z  
- ([y,z] \wedge x + [z,x]\wedge y)\Big) \in B_2(\g). $$ 

Next we note that $I_A$ is spanned by elements of the form 
$a \vee b - ab \vee \1,$ 
because $a \vee b \mapsto a\vee b - ab \vee \1$ is the projection of 
$S^2(A)$ onto $I_A$ with kernel $A \cong A \vee \1$. Therefore 
$B_2(\g)$ contains $I_A \otimes \k' \wedge \k$. On the other hand, 
(\ref{eq:sharp}) shows that $p_3(B_2(\g))$ 
is clearly contained in $I_A \otimes \k' \wedge \k$.
\end{Proof}

\begin{theorem} \label{theo1}
With the linear map 
$$ F \: A \otimes (\k \otimes \k \otimes \k) \to \Lambda^2(\g), 
\quad a \otimes (x\otimes y \otimes z) \mapsto (a \wedge \1 \otimes [x,y] \vee z) + 
a \otimes \partial(x \wedge y \wedge z) $$
we get the following description of $B_2(\g)$:  
$$ B_2(\g) = \Lambda^2(A) \otimes \k.S^2(\k) + T_0(A) \otimes \k \vee \k' 
+ \im(F) + I_A \otimes (\k \wedge \k').$$
\end{theorem}

\begin{proof} 
The description of the position of $B_2(\g)$ given in Lemma~\ref{lem2.3} is already 
quite detailed. It shows in particular that 
$$ B_2(\g) 
= (p_1 + p_2)(B_2(\g)) \oplus p_3(B_2(\g))
= (p_1 + p_2)(B_2(\g)) \oplus I_A \otimes (\k \cap \k') $$
and that $(p_1 + p_2)(B_2(\g))$ contains $\Lambda^2(A) \otimes \k.S^2(\k)$. 

We know from the proof of Lemma~\ref{lem2.3}(1) that, 
modulo the subspace $\Lambda^2(A) \otimes \k.S^2(\k) \subeq B_2(\g)$, we have 
\begin{eqnarray*}
(p_1 + p_2)(a x \wedge by \wedge cz)
&\equiv& T(a,b,c) \otimes [x,y] \vee z  + abc \otimes \partial(x \wedge y \wedge z)\\ 
&=& T_0(a,b,c) \otimes [x,y] \vee z  + abc \wedge \1 \otimes [x,y] \vee z + 
abc \otimes \partial(x \wedge y \wedge z) \\
&=& T_0(a,b,c) \otimes [x,y] \vee z  + F(abc \otimes x \otimes y \otimes z) 
\subeq T_0(A) \otimes \k \vee \k' + \im(F). 
\end{eqnarray*}
Since $T_0(A) \otimes \k \vee \k' \subeq B_2(\g)$ by Lemma~\ref{lem2.3}, 
we also obtain the converse inclusion 
$$ \im(F) \subeq (p_1 + p_2)(B_2(\g)) + T_0(A) \otimes \k \vee \k' + \Lambda^2(A) \otimes \k.S^2(\k) 
\subeq B_2(\g). $$
Now the theorem follows. 
\end{proof}

%%%%%%%%%%%%%%%%%%% new section %%%%%%%%%%%%%%%%%%%%%%%%%%%%%%
\section{The description of the $2$-cocycles} \label{sec3}

As explained in the introduction, 
elements of $H^2(\g,\z)$ can be identified with linear maps 
$f \: Z_2(\g) \to \z$, vanishing on the subspace $B_2(\g)$.
We further write $2$-cocycles as $f = f_1 + f_2 + f_3$, 
according to the decomposition in Lemma~\ref{lem2.1}, where 
$$ f_1 \: \Lambda^2(A) \otimes S^2(\k) \to \z, \quad 
f_2 \: A \otimes \Lambda^2(\k) \to \z \quad  \mbox{ and } \quad 
f_3 \: I_A \otimes \Lambda^2(\k) \to \z. $$ 
Here $f_1$ corresponds to an alternating bilinear map 
$\tilde f_1 \: A \times A \to \Sym^2(\k,\z)$, 
$f_2$ to a linear map 
$$\tilde f_2 \: A \to \Lin(\Lambda^2(\k),\z), \quad 
\tilde f_2(a)(x \wedge y) = \frac{1}{2}(f(ax \wedge y) - f(ay \wedge x)),$$  
and $f_3$ to a symmetric bilinear map 
$\tilde f_3 \: I_A \to \Alt^2(\k,\z).$
The condition, that three such maps $\tilde f_1, \tilde f_2, \tilde f_3$ combine 
to a $2$-cocycle 
$$f \: \Lambda^2(\g) \to \z, \quad 
a \wedge a' \otimes x \vee x' +  (b \vee \1) \otimes y  \wedge y'
+ c \otimes (z \wedge z') 
\mapsto \tilde f_1(a,a')(x,x') + \tilde f_2(b)(y,y') + \tilde f_3(c)(z,z'), $$
is that $f$ vanishes on $B_2(\g)$. To make this condition more explicit, 
we define the {\it Koszul map} 
$$ \Gamma \: \Sym^2(\k,\z)^\k  \to Z^3(\k,\z), \quad 
\Gamma(\kappa)(x,y,z) := \kappa([x,y],z). $$
That $\Gamma(\kappa)$ is alternating follows from 
$$ \Gamma(\kappa)(x,z,y) 
= \kappa([x,z],y) = \kappa(y,[x,z]) = \kappa([y,x],z) 
= - \Gamma(\kappa)(x,y,z) $$
and the fact that the symmetric group $S_3$ is generated by the transpositions 
$(1\ 2)$ and $(2\ 3)$. That the image of $\Gamma$ consists of $3$-cocycles 
is well known (\cite{Kos50}, \S 11; \cite{ChE48}, p.113).  

Recall that for each $\k$-module $\fa$, the Lie algebra differential 
$d_\k \: C^p(\k,\fa) \to C^{p+1}(\k,\fa)$
is given by 
\begin{eqnarray}
(d_\k \omega)(x_0, \ldots, x_p) 
&&:= \sum_{j = 0}^p (-1)^j x_j.\omega(x_0, \ldots, \hat x_j, \ldots, 
x_p) \cr
&& + \sum_{i < j} (-1)^{i + j} \omega([x_i, x_j], x_0, \ldots, \hat
x_i, \ldots, \hat x_j, \ldots, x_p), 
\end{eqnarray}
where $\hat x_j$ indicates omission of $x_j$. 

For the following theorem, we observe that the Lie algebra differential 
$d_\k \: C^2(\k,\z) = \Alt^2(\k,\z) \to Z^3(\k,\z)$ factors through the surjective map 
$\Alt^2(\k,\z) \onto \Lin(Z_2(\k),\z),$ whose kernel are the $2$-coboundaries. 

\begin{theorem}{\rm(Description of cocycles)}\label{theo2}
The function $f = f_1 + f_2 + f_3$ as above is a $2$-cocycle if and only if 
the following conditions are satisfied: 
\begin{description}
\item[(a)] $\im(\tilde f_1) \subeq \Sym^2(\k,\z)^\k$.  
\item[(b)] $\tilde f_1(T_0(A))$ vanishes on $\k \times \k'$. 
\item[(c)] $d_\k(\tilde f_2(a)) = \Gamma(\tilde f_1(a,\1))$ for each $a \in A$. 
\item[(d)] $\tilde f_3(I_A)$ vanishes on $\k \times \k'$. 
\end{description}
\end{theorem}

\begin{proof} The linear map $f$ is a $2$-cocycle if and only if it 
vanishes on $B_2(\g)$. In view of Theorem~\ref{theo1}, $B_2(\g)$ is 
the sum of four subspaces, so that we get four conditions. 

Condition (a) means that $f$ vanishes on $\Lambda^2(A) \otimes \k.S^2(\k)$, 
and condition (b) that it vanishes on the subspace $T_0(A) \otimes \k \vee \k'$. 

That $f$ vanishes on the image of $F$, means that 
$$ \Gamma(\tilde f_1(a,\1))(x,y,z) 
= \tilde f_1(a,\1)([x,y],z) 
= - \tilde f_2(a)(\partial(x \wedge y \wedge z)) 
= (d_\k\tilde f_2(a))(x,y,z) $$
for $a \in A$ and $x,y,z \in \k$, which is (c). 

Finally, (d) means that $f$ vanishes on $I_A \otimes \k \wedge \k'$.  
\end{proof}

\begin{corollary} $f = f_1 + f_2 + f_3$ is a cocycle if and only if 
$f_1 + f_2$ and $f_3$ are cocycles. 
\end{corollary} 

\begin{corollary} \label{cor2theo2} A function of one of the three types 
$f = f_i$, $i =1,2,3$,  
is a $2$-cocycle if and only if the following conditions are satisfied: 
\begin{description}
\item[($i=1$)] $\im(\tilde f_1) \subeq \Sym^2(\k,\z)^\k$ and the induced map 
$A \times A \to \Lin(\k \vee \k', \z)^\k$ is a cyclic $1$-cocycle. 
\item[($i=2$)] $\tilde f_2(A) \subeq Z^2(\k,\z)$. 
\item[($i=3$)] $\tilde f_3(I_A)$ vanishes on $\k \times \k'$. 
\end{description}
\end{corollary}

\begin{proof} That $f = f_i$ is a $2$-cocycle is equivalent to 
$f$ vanishing on $p_i(B_2(\g))$, so that Lemma~\ref{lem2.3} 
leads to the stated characterizations. 
\end{proof}

\begin{remark} \label{rem3.4} A special class of cocycles are those of the form 
$f = f_1$, vanishing on $\g \times \g'$. 
The cocycles of the form $f = f_3$ also vanish on the commutator algebra, and the 
sums of these two types exhaust the image of the injective pull-back map 
$H^2(\g/\g',\z)_{1,3} \cong \Alt^2(\g/\g',\z)_{1,3} \to H^2(\g,\z),$
where $\Alt^2(\g/\g',\z)_{1,3}$ denotes the set of all alternating 
maps vanishing on $(A \vee \1) \otimes (\k/\k' \wedge \k/\k') \subeq \Lambda^2(\g/\g')$. 
\end{remark}

\begin{corollary} \label{cortheo2} For each cocycle $f = f_1 + f_2 + f_3$,  
there exists a decomposition $f_1 = f_1^0 + f_1^1$, where 
$$ f_1^0(\g,\g') = \{0\}, \quad \im(\tilde f_1^1) \subeq \Sym^2(\k,\z)^\k 
\quad \mbox{ and } \quad T_0(A) \subeq \ker \tilde f_1^1. $$
\end{corollary}

\begin{proof} Conditions (a) and (b) in Theorem~\ref{theo2} only refer to the restriction $\oline f_1^1$ of 
$f_1$ to the subspace $\Lambda^2(A) \otimes (\k \vee \k')$ of $\Lambda^2(A) \otimes S^2(\k)$. 
This has the following interesting consequence. We have a short exact sequence 
$$ \{0\} \to \Sym^2(\k/\k',\z) \to \Sym^2(\k,\z)^\k \to \Lin(\k \vee \k',\z)^\k \to \{0\}, $$
where the surjectivity of the map $\Sym^2(\k,\z)^\k \to \Lin(\k \vee \k',\z)^\k$ 
follows from the fact that any symmetric bilinear extension of an element of 
$\Lin(\k \vee \k',\z)^\k$ is invariant. 
Any splitting of this sequence extends $\oline f^1_1$ to an alternating bilinear map 
$\tilde f_1^1 \: A \times A \to \Sym^2(\k,\z)^\k$
with 
$$ \tilde f_1^1(a,b)(x,y) = \tilde f_1(a,b)(x,y) 
\quad \mbox{ for } \quad a,b \in A, x \in\k,y \in \k'$$
and such that $T_0(A) \subeq \ker\tilde f_1^1$. Then 
$$ \Gamma(\tilde f_1^1(a,\1)) =  \Gamma(\tilde f_1(a,\1)) \quad \mbox{ for } \quad a \in A, $$
so that $f_1^1 + f_2 + f_3$ also is a cocycle by Theorem~3.1. 
We conclude that 
$f_1^0 := f_1 - f_1^1$ 
is a cocycle vanishing on $\g \times \g'$. 
This proves the assertion. 
\end{proof}

\begin{proposition}{\rm(Description of coboundaries)} \label{prop-cobound}
A cocycle $f = f_1 + f_2 + f_3$ is a coboundary if and only if 
$f_1 = f_3 = 0$ 
and there exists a linear map $\ell \: A \to \Lin(\k,\z)$ with 
$d_\g \ell = f_2$, i.e., 
$$ \tilde f_2(a) = d_\k(\ell(a)) \quad \mbox{ for all } \quad a \in A. $$
\end{proposition}

\begin{proof} That $f$ is a coboundary means that it vanishes on 
$Z_2(\g)$. According to Lemma~\ref{lem2.1}, this implies that 
$f_1 = f_3 = 0$. Since the bracket map $b_\g \: \Lambda^2(\g) \to \g$ 
is alternating in $\k$ and symmetric in $A$, all coboundaries are of the form 
$f = f_2$.   
\end{proof}

A {\it coupled cocycle} is a cocycle of the form $f_1 + f_2$ for which $f_1$ is not a cocycle. 
The following theorem characterizes the pairs 
$(A,\k)$ for which $A \otimes \k$ possesses coupled cocycles. 
In Section~\ref{sec5} below, we shall also give a concrete example of a Lie algebra 
$\k$ satisfying this condition. 

\begin{theorem} \label{theo3} The Lie algebra $\g = A \otimes \k$ possesses coupled cocycles 
if and only if $d_A(A) \not=\{0\}$ and $\k$ possesses a symmetric invariant bilinear 
form $\kappa$ for which $\Gamma(\kappa) \in Z^3(\k)$ is a non-zero coboundary. 

If this is not the case, then each cocycle $f \in Z^2(\g)$ is a sum 
$$f = f_1 + f_2 + f_3 = f_1^0 + f_1^1 + f_2 + f_3. $$
of four cocycles, where 
\begin{description}
\item[(a)] $f_1^0$ vanishes on $\g\times \g'$. 
\item[(b)] $\tilde f_1^1 \in Z^1(A, \Sym^2(\k,\z)^\k)$ is a cyclic $1$-cocycle.  
\item[(c)] $\tilde f_2(A) \subeq Z^2(\k,\z)$. 
\item[(d)] $f_3$ vanishes on $\g \times \g'$. 
\end{description}
\end{theorem}

\begin{proof} First let $f = f_1 + f_2$ be a coupled cocycle on $\g$. Then we have 
$\Gamma(\tilde f_1(A,\1)) \not= \{0\}$. 
Composing with a suitable linear functional $\chi \: \z \to \K$ 
with 
$$ \chi \circ \Gamma(\tilde f_1(A,\1)) 
= \Gamma((\chi \circ f_1)\,\tilde{}\,\,(A,\1)) \not= \{0\}, $$
we may w.l.o.g.\ assume that $\z = \K$. 
Then there exists an $a \in A$ with 
$$ d_\k(\tilde f_2(a)) = \Gamma(\tilde f_1(a,\1))\not=0. $$ 
Now $\kappa := \tilde f_1(a,\1) \in \Sym^2(\k)^\k$ is an invariant symmetric 
bilinear form for which $\Gamma(\kappa)$ is exact and non-zero. 
Then $a \wedge \1 \in T(A) \setminus T_0(A)$ (Theorem~\ref{theo2}), 
so that $0 \not= d_A(a)$ in $\Omega^1(A)$ (Lemma~\ref{lem1.1}). 

If, conversely, $d_A(A) \not=\{0\}$ and $\kappa$ is an invariant symmetric 
bilinear form on $\k$ for which $\Gamma(\kappa)$ is a non-zero coboundary, 
then we pick $\eta \in C^2(\k) = \Alt^2(\k)$ with $d_\k \eta = \Gamma(\kappa)$. 
We now define linear maps 
$$ \tilde f_1 := \gamma_A \otimes \kappa 
\: \Lambda^2(A) \to \Sym^2(\k,\Omega^1(A))^\k, \quad 
\tilde f_1(a \wedge b)(x,y) := \kappa(x,y)\cdot (a d_A(b)- b d_A(a)) $$
and 
$$ \tilde f_2 := -d_A \otimes \eta 
\: A \to C^2(\k,\Omega^1(A)), \quad 
\tilde f_2(a)(x,y) := -\eta(x,y) \cdot d_A(a). $$

We claim that the corresponding map $f = f_1 + f_2$ is a $2$-cocycle 
by verifying the conditions in Theorem~\ref{theo2}. 
Condition (a) is obviously satisfied, and (b) follows from 
$T_0(A) = \ker \gamma_A$ (Lemma \ref{lem1.1}). Further $f_3 = 0$, and (c) follows from 
$$ d_\k \tilde f_2(a) = -(d_\k \eta) \cdot d_A(a) 
= -\Gamma(\kappa) d_A(a) = \Gamma(\tilde f_1(a,\1)). $$

That $f_1$ is  not a cocycle, i.e., that $f$ is coupled, means that 
$\tilde f_1(A \wedge \1)(\k \times \k') \not=\{0\}$, which is equivalent to 
$d_A(A) \not= \{0\}$  and $\Gamma(\kappa) = \eta \not=0.$
This completes the proof of the first part of the theorem. 

For the second part, we assume that either $d_A(A) \cong T_0(A)/T(A)$ 
vanishes, which means that $T_0(A) = T(A)$, or that 
for each exact invariant symmetric bilinear form $\kappa$ on $\k$ we have $\Gamma(\kappa) = 0$. 
Then for each cocycle $f = f_1^0 + f_1^1 + f_2 +f_3$ as in 
Corollary~\ref{cortheo2}, either 
$\tilde f_1^1$ vanishes on $T(A)$ (if $d_A(A)$ vanishes) 
or $\tilde f_2(A) \subeq Z^2(\k,\z)$ (if for all exact forms on $\k$ 
the $3$-cocycle $\Gamma(\kappa)$ vanishes). Both 
conditions imply that $f^1_1$ and $f_2$ are cocycles. 
Hence the assertion follows from Corollary~\ref{cor2theo2}. 
\end{proof} 

\begin{corollary} If $H^1(\k,\k^*)= \{0\}$, then $\g = A \otimes \k$ has no coupled cocycles. 
\end{corollary}

\begin{proof} From the exact sequence in Proposition~\ref{transfer-seq} below, 
it follows that 
the Koszul map 
$$\gamma \: \Sym^2(\k)^\k \to H^3(\k), \kappa \mapsto [\Gamma(\kappa)]$$ 
is injective, and this implies that 
each exact invariant form vanishes.   
\end{proof}

The following proposition describes the universal cocycle for $\g$ in terms of our 
threefold direct sum decomposition. 

\begin{proposition}{\rm(A universal cocycle)} 
Let $p_\k \: \Lambda^2(\k) \to Z_2(\k)$ denote a linear projection onto 
$Z_2(\k)$. Then the linear map 
$$ \tilde f^u := p_1 \oplus (\id_A \otimes p_\k) \oplus p_3 \: \Lambda^2(\g) \to 
Z_2(\g) $$
maps $B_2(\g)$ into itself, hence induces a $2$-cocycle 
$$ f^u \: \Lambda^2(\g) \to H_2(\g) = Z_2(\g)/B_2(\g). $$
It is universal in the sense that for each space $\z$ the map 
$$ \Lin(H_2(\g),\z) \to H^2(\g,\z), \quad \phi \mapsto \phi  \circ f^u $$
is a linear bijection. 
\end{proposition}

\begin{proof} That $\tilde f^u$ is a linear projection onto $Z_2(\g)$ 
follows from Lemma~\ref{lem2.1}. The remainder follows from the fact that 
$H^2(\g,\z) \to \Lin(Z_2(\g),\z),  [f] \mapsto f\res_{Z_2(\g)}$
is injective onto the set of all maps vanishing on $B_2(\g)$. 
\end{proof}

%%%%%%%%%%%%%%%%%%%%%%%%%%%%%%%%%%%%%%%%%%%%%%%%%%%%%%%%%%%%%%%%%%%%%%
\section{The structure of the second cohomology space} \label{sec4} 

In this section, we use the results of the previous section to give a quite explicit description of 
the space $H^2(\g) = H^2(\g,\K)$ 
in terms of data associated directly to $\g$ and $A$. 

\begin{lemma} \label{lem4.1} Associating with each linear map 
$\tilde f_2 \: A \to Z^2(\k)$ the corresponding cocycle $f_2 \in Z^2(\g)$, we obtain, 
together with the natural pull-back map $H^2(\g/\g') \to H^2(\g)$, an injection 
$$ H^2(\g/\g')_{1,3} \oplus \Lin(A, H^2(\k)) \ssmapright{\Phi} H^2(\g) $$
whose image consists of all classes of cocycles of the form 
$f_1^0 + f_2 + f_3$.   
\end{lemma}

\begin{proof} The image of the pull-back map 
$H^2(\g/\g')_{1,3} \to H^2(\g)$ consists of those 
cohomology classes represented by cocycles vanishing on $\g \times \g'$, 
which are the cocycles of the form $f_1^0 + f_3$. Since the space of these cocycles
intersects $B^2(\g)$ trivially, the space $H^2(\g/\g')_{1,3}$ injects into $H^2(\g)$ 
(Remark~\ref{rem3.4} and Prop.~\ref{prop-cobound}). 

Next we recall that the cocycles of the form $f = f_2 \: A \otimes \Lambda^2(\k) \to \K$ 
correspond to linear maps $\tilde f_2 \: A \to Z^2(\k)$ (which means that 
$f_2$ vanishes on $A \otimes B_2(\k)$), and that such a map is a 
coboundary if and only if $\im(\tilde f_2)(A) \subeq B^2(\k)$, because this implies 
the existence of a linear map $\ell \: A \to \Lin(\k)$ with 
$\tilde f_2(a) = d_\k(\ell(a))$ for all $a \in A$. The latter condition means that 
$f_2$ vanishes on $A \otimes Z_2(\k)$, so that the cohomology classes correspond to 
elements in 
$$ \Lin(A \otimes Z_2(\k)/(A \otimes B_2(\k)),\K) 
\cong \Lin(A \otimes (Z_2(\k)/B_2(\k)),\K)
\cong \Lin(A \otimes H_2(\k), \K)
\cong \Lin(A, H^2(\k)). $$
\end{proof}

Given a cocycle $f = f_1 + f_2 + f_3$ in $Z^2(\g)$, we obtain the map 
$\Gamma \circ \tilde f_1 \: \Lambda^2(A) \to Z^3(\k),$
whose kernel contains $T_0(A)$, so that it induces a linear map 
$$ f^\flat \: \Omega^1(A) \cong \Lambda^2(A)/T_0(A) \to Z^3(\k), \quad 
a \cdot d_A(b) - b \cdot d_A(a) \mapsto \Gamma(\tilde f_1(a,b)), $$
mapping the subspace $d_A(A) \subeq \Omega^1(A)$ into the subspace 
$B^3(\k)$ (Theorem~\ref{theo2}). In view of 
\begin{eqnarray}\label{eq1}
f^\flat(d_A(a)) = -\Gamma(\tilde f_1(a,\1)) = -d_\k(\tilde f_2(a)), 
\end{eqnarray}
the range of each map $\Gamma \circ \tilde f_1$ 
lies in the subspace $Z^3(\k)_\Gamma := \im(\Gamma) \subeq Z^3(\k)$ and 
$$ f^\flat(d_A(A)) \subeq B^3(\k)_\Gamma := B^3(\k) \cap \im(\Gamma). $$
We thus obtain a map 
$$ \Psi \: H^2(\g) \to \Lin((\Omega^1(A), d_A(A)), (Z^3(\k)_\Gamma, B^3(\k)_\Gamma)), 
\quad [f] \mapsto \Gamma \circ \tilde f_1, $$
where for pairs $(X,X')$ and $(Y,Y')$ of linear spaces with $X' \subeq X$ and 
$Y' \subeq Y$ we write 
$$ \Lin((X,X'), (Y,Y')) := \{ f \in \Lin(X,Y) \: f(X') \subeq Y'\}. $$

\begin{theorem} \label{theo4} The sequence 
$$ \{0\} \to H^2(\g/\g')_{1,3} \oplus \Lin(A, H^2(\k)) \ssmapright{\Phi} H^2(\g) 
\ssmapright{\Psi} \Lin((\Omega^1(A), d_A(A)), (Z^3(\k)_\Gamma, B^3(\k)_\Gamma)) \to \{0\}$$ 
is exact. 
\end{theorem}

\begin{proof} We have already seen in Lemma~\ref{lem4.1} that $\Phi$ is injective. 

The kernel of $\Psi$ consists of all cocycles $f = f_1 + f_2 + f_3$ for which 
$\Gamma \circ \tilde f_1 = 0$. This is equivalent to 
$\tilde f_1(\Lambda^2(A))$ vanishing on $\k \vee \k'$, which means that 
$f_1$ vanishes on $\g \times \g'$, i.e., $f_1 = f_1^0$. We conclude that  
$\ker \Psi = \im \Phi$. 

To see that $\Psi$ is surjective, let 
$\alpha \in \Lin((\Omega^1(A), d_A(A)), (Z^3(\k)_\Gamma, B^3(\k)_\Gamma))$ 
and observe that there exists a linear map 
$$ f^\flat \: \Omega^1(A) \to \Sym^2(\k)^\k 
\quad \mbox{ with } \quad 
\Gamma \circ f^\flat = \alpha, $$
and a linear map $\beta \: d_A(A) \to C^2(\k)$ with 
$$ d_\k(\beta(d_A(a))) = \alpha(d_A(a)) \quad \mbox{ for all } \quad a \in A.$$ 
For 
$$ \tilde f_1 \: \Lambda^2(A) \to \Sym^2(\k)^\k,\quad 
\tilde f_1(a,b) := f^\flat(a \cdot d_A(b)-b \cdot d_A(a))
\quad \mbox{ and } \quad 
\tilde f_2 \: A \to C^2(\k), \quad a \mapsto -\beta(d_A(a)),$$
we then have 
$$ d_\k(\tilde f_2(a)) = -d_\k(\beta(d_A(a))) =- \alpha(d_A(a)) 
= -\Gamma(f^\flat(d_A(a))) 
= \Gamma(\tilde f_1(a,\1)), $$
so that the corresponding maps $f_1$ and $f_2$ sum up to a $2$-cocycle 
$f := f_1 + f_2$ satisfying $\Psi([f]) = \alpha$. 
\end{proof}

The quotient
$Z^3(\k)_\Gamma/B^3(\k)_\Gamma$ can be identified with the image of the map 
$$ \gamma \: \Sym^2(\k)^\k \to H^3(\k), \quad \kappa \mapsto [\Gamma(\kappa)] $$
discussed in the appendix below. From the exactness of the sequence 
in Proposition~\ref{transfer-seq}, it follows that the space 
$\Sym^2(\k)^\k_{\rm ex} := \ker\gamma$ of exact invariant bilinear forms 
satisfies 
\begin{eqnarray}
  \label{eq:sym}
\Sym^2(\k)^\k_{\rm ex} \cong  H^1(\k,\k^*)/H^2(\k). 
\end{eqnarray}

We also note that for a {\it quadratic Lie algebra}, i.e., a finite-dimensional 
Lie algebra $\k$ with an invariant non-degenerate symmetric bilinear form $\kappa_0$, the space 
$\out(\k) := \der(\k)/\ad \k$ of {\it outer derivations} satisfies 
$$ H^1(\k,\k^*) \cong H^1(\k,\k) \cong \der(\k)/\ad \k = \out(\k), $$
and that the subspace $H^2(\k) \subeq H^1(\k,\k^*)$ consists of those classes $[D]$ 
of derivations $D$ which are skew-symmetric with respect to $\kappa_0$. 

We further have $\ker \Gamma \cong \Sym^2(\k/\k')$, so that 
$$ B^3(\k)_\Gamma \cong \Sym^2(\k)^\k_{\rm ex}/\Sym^2(\k/\k') 
\quad \mbox{ and } \quad 
Z^3(\k)_\Gamma \cong \Sym^2(\k)^\k/\Sym^2(\k/\k'). $$
To obtain an explicit description of $H^2(\g)$, it is therefore necessary to have a good 
description of the space $\Sym^2(\k)^\k$ of invariant quadratic forms on $\k$ and its 
subspace of exact forms.

\begin{problem} Let $\k$ be a finite-dimensional $\K$-Lie algebra. 
We consider the space $S := \Sym^2(\k)^\k$ of invariant symmetric bilinear forms on $\k$. 

Let $\fn := \bigcap \{ \rad(\kappa) \: \kappa \in S\}$ denote the common radical of all 
invariant symmetric bilinear forms on $\k$. 
Fix an element $\kappa \in S$ of maximal rank. Then 
$\fn \subeq \rad(\kappa)$, but is there some $\kappa$ for which we 
have equality? 
\end{problem}

In the following remark, we collect some information that is useful to determine 
the space $Z^3(\k)_\Gamma$. 
\begin{remark} Suppose that $(\k,\kappa_0)$ is a quadratic Lie algebra, i.e., 
$\kappa_0$ is a non-degenerate invariant symmetric bilinear form on $\k$. 
Then there exists for each invariant symmetric bilinear form 
$\kappa \in \Sym(\k)^\k$ a uniquely determined endomorphism 
$A_\kappa \in \End(\k)$ with 
$$ \kappa(x,y) = \kappa_0(A_\kappa.x, y) \quad \mbox{ for } \quad x,y \in \k. $$
Now the invariance of $\kappa$ implies that 
$A_\kappa$ is contained in the {\it centroid} 
$$ \Cent(\k) := \{A \in \End(\k) \: (\forall x \in \k)\ [A, \ad x] = 0\}. $$
The centroid of $\k$ is an associative subalgebra of $\End(\k)$ 
on which transposition $A \mapsto A^\top$ with respect to $\kappa_0$ induces 
a linear anti-automorphism, satisfying 
$$ \kappa_0(A.x,y) = \kappa_0(x, A^\top.y) = \kappa_0(A^\top.y, x)\quad \mbox{ for all} 
\quad x,y \in \k. $$
It follows in particular that for $A \in \Cent(\k)$ the invariant bilinear form 
$\kappa_A(x,y) := \kappa_0(A.x, y)$ is symmetric if and only if $A^\top = A$. 
This leads to a linear bijection 
$$ \Cent(\k)_+ := \{ A \in \Cent(\k) \: A^\top = A\} \to \Sym^2(\k)^\k, \quad 
A \mapsto \kappa_A. $$

For $A^\top = -A$, the invariant form $\kappa_A$ is alternating, 
which implies that $\kappa_A$ vanishes on $\k \times \k'$, and this implies that 
$$A(\k) \subeq (\k')^\bot = \z(\k) 
\quad \mbox{ and } \quad A(\k') = \{0\}. $$
Conversely, any $A \in \End(\k)$ with $\k' \subeq \ker A$ and 
$\im(A) \subeq \z(\k)$ satisfies $A \circ \ad x = \ad x \circ A = 0$ for 
all $x \in \k$, hence is contained in the centroid. We put 
$$ \Cent_0(\k) := \{ A \in \End(\k) \: \k' \subeq \ker A, \ \im(A) \subeq \z(\k)\} $$
and observe that 
$\Cent_0(\k) \trile \Cent(\k)$
is an ideal of the associative algebra $\Cent(\k)$ because 
$$ \Cent_0(\k) = \{ A \in \Cent(\k) \: A\res_{\k'} = 0\} $$
is the kernel of the restriction homomorphism $\Cent(\k) \to \End(\k')$. 

If $A \in \Cent_0(\k)$, then 
$$ \kappa_0(A^\top.[\k,\k],\k) = \kappa_0([\k,\k],A.\k) 
\subeq \kappa_0(\k', \z(\k)) = \{0\}, $$
so that $A^\top \in \Cent_0(\k)$. Hence the ideal 
$\Cent_0(\k)$ is invariant under transposition. We have already seen that 
$\Cent_0(\k)$ contains all skew-symmetric elements of $\Cent(\k)$, so that the 
involution induced on the quotient algebra 
$$ \Cent_{\rm red}(\k) := \Cent(\k)/\Cent_0(\k) \into \End(\k') $$ 
is trivial, which implies that this algebra is commutative. 

We thus have 
$$ \Cent(\k)_+ \cong \Sym^2(\k)^\k \quad \mbox{ and } \quad 
\Cent_0(\k)_+ := \{ A \in \Cent_0(\k) \: A^\top = A \} \cong \Sym^2(\k/\k'). $$
Therefore 
$$ Z^3(\k)_\Gamma = \im(\Gamma) 
\cong \Cent(\k)_+/\Cent_0(\k)_+ 
\cong \Cent_{\rm red}(\k) $$
carries the structure of an associative commutative algebra. 

In \cite{MR93}, Th.~2.3, Medina and Revoy describe the structure 
of the associative algebra 
$\Cent(\fk)$ for a Lie algebra $\fk$ whose center $Z(\fk)$ is contained in 
$\fk'$: The algebra $\Cent(\fk)$ has a decomposition with respect to 
orthogonal indecomposable idempotents $e_1, \ldots, e_r$ with 
$\sum_i e_i = \id_\k$, so that $\fk$ is the direct product of the ideals $\fk_i :=e_i\fk$. 
Moreover, the algebra $\Cent(\fk_i)\simeq e_i\Cent(\fk) e_i$ is a local ring, and we have 
$$ \Cent(\fk)=\bigoplus_{i,j=1}^r \Cent_{ij}, \quad \mbox{ where } \quad \Cent_{ij} := e_i\Cent(\fk) e_j \cong \Lin(\fk_j/\fk_j',Z(\fk_i)), i \not=j, $$ 
as linear spaces, and 
$$\Cent_0(\fk)=\Big(\bigoplus_{i=1}^r \Cent_0(\fk_i)\Big)
\oplus\Big(\bigoplus_{i\not= j}\Cent_{ij}\Big). $$
If, in addition, $\k$ carries a non-degenerate quadratic from 
$\kappa_0$, then Th. 2.5 {\it loc.cit.} implies that 
the decomposition of $\fk$ as a direct sum of ideals $\fk_i$ is 
orthogonal and the idempotents $e_i$ are symmetric with respect to $\kappa_0$. 
We conclude in particular that 
$$ \Cent_{\rm red}(\k) = \Cent(\k)/\Cent_0(\k)\cong \bigoplus_{i=1}^r \Cent_{\rm red}(\k_i). $$
\end{remark}

%%%%%%%%%%%%%%%%%%% new section %%%%%%%%%%%%%%%%%%%%%%%%%%%%%%
\section{Some examples} \label{sec5}

In this section we describe some Lie algebras $\k$ on which we have 
invariant bilinear forms $\kappa$ for which $\Gamma(\kappa)$ is a non-zero coboundary, 
so that $\g = A \otimes \k$ has coupled cocycles whenever $d_A\not=0$. 

\subsection{The split oscillator algebra} 

Let $\h$ be the $3$-dimensional Heisenberg algebra 
$\fh$ with generators $x$, $y$ and $c$ and the only non-trivial relation $[x,y]=c$.
 Then pass to the extension $\k = \fh \rtimes \K D$ of $\fh$ 
by a derivation $D$ like for affine Kac--Moody algebras. 
Explicitly, we take $D(x)=x$, $D(y)=-y$ and $D(c)=0$  
(cf. \cite{MooPia}, p.98, Ex.\ $6$). 
We write $d := (0,1)$ for the element of 
$\k$ corresponding to $D$. 
The Lie algebra $\k$ 
is $4$-dimensional, and has an invariant bilinear symmetric form $\kappa$, 
as any Lie algebra with symmetrizable Cartan matrix (cf. \cite{MooPia}, 
Prop.~4, p.\ 362).
We call $\k$ the {\it split oscillator algebra over $\K$}. 
\begin{remark}
Let us compute the dimensions of the spaces of cochains, cocycles and 
cohomology spaces:

\vspace{1cm}\hspace{4cm}
\begin{tabular}{|c|c|c|c|c|c|}  \hline
degree $p$ & 0 & 1 & 2 & 3 & 4 \\   \hline\hline
$\dim C^p(\k)$ & 1 & 4 & 6 & 4 & 1 \\ \hline
$\dim H^p(\k)$ & 1 & 1 & 0 & 1 & 1 \\ \hline
$\dim B^p(\k)$ & 0 & 0 & 3 & 3 & 0 \\ \hline
$\dim Z^p(\k)$ & 1 & 1 & 3 & 4 & 1 \\ \hline
\end{tabular}\vspace{.5cm}  

In the preceding table, the dimension of the cohomology spaces is computed as follows: 
$\dim H^0(\k) = 1$ by definition. 
As $\k/[\k,\k]\,=\,\K D$, $\dim H^1(\k)\,=\,1$. By unimodularity, $\k$ satisfies 
Poincar\'e duality (\cite{Fu86}, p.~27), so that 
the dimensions in degree $3$ and $4$ follow. But the Euler characteristic of 
a finite dimensional Lie algebra vanishes \cite{Gol}, which implies that 
$H^2(\k) = \{0\}$. 

The dimensions of the boundary spaces are clear in degree $0$ and $1$. In degree $2$, there remain $3$ dimensions as the difference of  $\dim C^1(\k)$ and  $\dim Z^1(\k)$. In the same way, we get 
the dimensions of $B^p(\k)$ for $p=3,4$. Finally, $\dim Z^p(\k)$ is the sum of $\dim B^p(\k)$ and $\dim H^p(\k)$.

Observe that $[\k,\k]=\fh$ and 
$[\fh,\fh]=\K c$, so that $\k$ is solvable, 
but $[\k,\fh]=\fh$, so that $\k$ is not nilpotent.

We claim that each invariant bilinear form $\kappa$ is exact, which gives rise to coupled 
cocycles (in the sense of Section $3$):
If $0 \not=\mu \in C^4(\k)$, then the fact that $\k$ is unimodular 
implies that all $3$-cochains $i_h \mu$, $h \in \k$, are $3$-cocycles. 
If $h \in [\k,\k] = \h$, then $i_h \mu$ is exact, so that 
$i_D \mu$ yields a basis of the one-dimensional space $H^3(\k)$. 
Since $0 \not= (i_D \mu)(x,y,c)= \mu(D,x,y,c)$ and for each invariant symmetric 
bilinear form $\kappa$ we have 
$\kappa([x,y],c)
=\kappa(x,[y,c])=0$, we see that $\Gamma(\kappa) \in 
\spann \{ i_h \mu \: h \in \h\} = B^3(\k)$. Hence each invariant symmetric 
bilinear form is exact. 
\end{remark}

\begin{remark} We now turn to the space $\Sym^2(\k)^\k$: 
Any invariant symmetric bilinear form $\kappa$ satisfies 
$$\kappa(c,x)\,=\,\kappa([x,y],x)\,=\,-\kappa([y,x],x)\,=\,-\kappa(y,[x,x])\,=\,0.$$ 
$$\kappa(c,y)\,=\,\kappa([x,y],y)\,=\,\kappa(x,[y,y])\,=\,0.$$
$$\kappa(c,c)\,=\,\kappa([x,y],c)\,=\,\kappa(x,[y,c])\,=\,0.$$
$$\kappa(d,x)\,=\,\kappa(d,[d,x])\,=\,\kappa([d,d],x)\,=\,0.$$
$$\kappa(d,y)\,=\,-\kappa(d,[d,y])\,=\,-\kappa([d,d],y)\,=\,0.$$
$$\kappa(d,c)\,=\,\kappa(d,[x,y])\,=\,\kappa([d,x],y)\,=\,\kappa(x,y).$$
$$\kappa(x,x)\,=\,\kappa(x,[d,x])\,=\,-\kappa(x,[x,d])\,=\,-\kappa([x,x],d)\,=\,0.$$
$$\kappa(y,y)\,=\,-\kappa(y,[d,y])\,=\,\kappa(y,[y,d])\,=\,\kappa([y,y],d)\,=\,0.$$
We immediately conclude that the space of invariant symmetric bilinear forms 
is at most $2$-dimensional and that each such form $\kappa$ is determined by  
$\kappa(d,c)=\kappa(x,y)$ and $\kappa(d,d)$ (note that $d$ is not a commutator). 
Let us denote by $\kappa_1$ the 
(invariant symmetric bilinear form) with $\kappa_1(d,d)=1$ and $\kappa_1(x,y) = 0$ 
and $\kappa_2$ the 
invariant symmetric bilinear form with $\kappa_2(x,y)=\kappa_2(d,c)=1$ and 
$\kappa_2(d,d) = 0$. Then $\kappa_2$ coincides with the invariant form 
$\kappa$ introduced above and $\kappa_1, \kappa_2$ form a basis of $\Sym^2(\k)^\k$. 
Combining with the observation in the preceding remark and Section 4, we get 
$$ \K^2 \cong \Sym^2(\k)^\k \cong \Sym^2(\k)^\k_{\rm ex} \cong H^1(\k,\k^*) 
\cong H^1(\k,\k) \cong \out(\k). $$

For the reduced centroid, we thus get 
$${\rm Cent}_{\rm red}(\fk)\,=\,{\rm Sym}^2(\fk)^{\fk}\,
/\,{\rm Sym}^2(\fk\,/\,[\fk,\fk]) = \K[\kappa_2].$$
We further get $Z^3({\mathfrak k})_{\Gamma}=B^3({\mathfrak k})_{\Gamma} \cong \K$.

For any algebra $A$, and $\g = A \otimes \k$,  
the exact sequence in Theorem~4.2 now turns into a sequence of the form 
$$ \{0\} \to H^2(\g/\g') \cong \Lambda^2(A)^* \ssmapright{\Phi} 
H^2(\g) \to \Lin(\Omega^1(A), Z^3(\k)_\Gamma) \cong \Omega^1(A)^* \to \{0\}. $$
Therefore the essential part of $H^2(\g)$ is isomorphic to the dual space of 
$\Omega^1(A)$. From the construction in the proof of Theorem~4.2, it follows that 
the coupled cocycles correspond to the elements of $\Omega^1(A)^*$ not vanishing on 
the subspace $d_A(A)$. 
\end{remark}

\subsection{Two more classes of examples} 

\begin{example} For a Lie algebra $\g$ its cotangent bundle  
$T^*\g := \g^* \rtimes \g$ is a Lie algebra with the bracket 
$$ [(f,x), (f',x')] := (x.f' - x'.f,[x,x']), 
\quad \hbox{ where } \quad x.f = \ad^*x.f = -f \circ \ad x. $$
A slight generalization is obtained as follows. 
Suppose that $\gamma \: \g \times \g \to \g^*$ is a 
Lie algebra $2$-cocycle, i.e., 
$$ \sum_\cyc \ad^*x.\gamma(y,z) - \gamma([x,y],z) = 0. $$
Then we have another Lie algebra structure on $\g^* \oplus \g$ given by 
$$ [(f,x), (f',x')] := (x.f' - x'.f + \gamma(x,x'),[x,x']). $$
We write $T^*_\gamma \g$ for the corresponding Lie algebra, a so-called 
{\it twisted magnetic extension of $\g$}. 

The symmetric bilinear form given by 
$$ \kappa((f,x), (f',x')) := f(x') + f'(x) $$
satisfies 
\begin{eqnarray}
\Gamma(\kappa)((f,x), (f',x'), (f'',x'')) 
&=& \kappa((x.f' - x'.f + \gamma(x,x'),[x,x']), (f'',x'')) \cr
&=& f''([x,x']) +f'([x'',x]) + f([x',x'']) + \gamma(x,x')(x'') \cr
&=& \gamma(x,x')(x'') + \sum_\cyc f([x',x'']). 
\end{eqnarray}
This implies that $\kappa$ is invariant if and only if $\tilde\gamma(x,x',x'') := 
\gamma(x,x')(x'')$ is alternating, 
hence an element of $Z^3(\g)$ (cf.\ \cite{Bo97}, Lemma~3.1). 
Let us assume that this is the case and note that $\Gamma(\kappa)$ vanishes only 
if $\gamma=0$ and $\g$ is abelian. 

For the alternating bilinear form $\eta((f,x), (f',x')) := f(x') - f'(x)$
we then have 
\begin{eqnarray}
\eta([(f,x), (f',x')], (f'',x'')) 
&=& (x.f' - x'.f + \gamma(x,x'))(x'') - f''([x,x']) \cr
&=& \gamma(x,x')(x'') + f'([x'',x]) + f([x',x'']) - f''([x,x']), 
\end{eqnarray}
so that 
\begin{eqnarray}
(d \eta)((f,x), (f',x'), (f'',x'')) 
&=& -  \sum_\cyc\big(f'([x'',x]) + f([x',x'']) - f''([x,x'])\big) 
- \sum_\cyc \gamma(x,x')(x'') \cr
&=& - 3 \gamma(x,x')(x'') -  \sum_\cyc  f([x',x'']). 
\end{eqnarray}

Let $q \: T^*_\gamma \g \to \g$ denote the canonical projection. 
Then the preceding calculation shows that $\Gamma(\kappa)$ is a coboundary 
if and only if $[q^*\tilde\gamma] \in H^3(T^*_\gamma \g)$ vanishes. 
This is in particular the case for $\gamma = 0$. 
\end{example} 

\begin{example} In \cite{Pe97}, Pelc introduces a family of Lie algebras ${\cal A}_n$, 
where ${\cal A}_n$ is an $(n+1)$-dimensional Lie algebra 
with basis $T_0, \ldots, T_n$ and commutator relations 
$$ [T_i, T_j] = \left\{ 
  \begin{array}{cl} 
\hat{i-j} \cdot T_{i+j} & \mbox{for $i + j \leq n$} \\ 
0 &\mbox{ otherwise.}  
\end{array} \right. $$
Here $\hat i \in \{-1,0,1\}$ is chosen such that $i - \hat i \in 3 \Z$. 

Then ${\cal A}_n$ is a solvable Lie algebra with commutator algebra 
${\cal A}_n' = \spann \{ T_i \: i > 0\}$. Let us assume that $n = 3m$ for some $m \in \N$ and that 
$\ch\K = 0$. 
Then ${\cal A}_{3m}$ carries a non-degenerate invariant 
symmetric bilinear form defined by 
$$ \kappa(T_i, T_j) = \delta_{i+j,n} $$
(cf.\ \cite{Pe97}). 
We claim that $\kappa$ is exact. Choose $a_0, \ldots, a_n \in \K$ in such a way 
that $a_i = \alpha i + \beta$, where $\beta = -1$ and $\alpha = \frac{2}{n}$. 
Then $a_0 = -1$ and $a_i = -a_{n-i}$, and we define a $2$-cochain 
$\eta \in C^2({\cal A}_n,\K)$ by 
$$ \eta(T_i, T_j) = a_i \delta_{j,n-i}. $$
Note that we need $a_i = - a_{n-i}$ to see that this is well-defined. 

Now $\eta([T_i, T_j], T_k)$ vanishes if $i + j + k \not= n$, and 
for $i + j + k = n$ we get 
$$ \eta([T_i, T_j], T_k) = \hat{i-j} \cdot a_{i+j}. $$
With $\hat{j-k} = (j - (n - j-i))\,\hat{} = \hat{i -j},$ 
this leads to 
$$ -d \eta(T_i, T_j, T_k) 
= \hat{i-j} a_{i+j} + \hat{j-k} a_{j+k} + \hat{k-i} a_{k+i} 
= \hat{i-j}(a_{i+j} - a_i -  a_j)
= \hat{i-j} = \kappa([T_i, T_j], T_k). $$
Therefore $\kappa$ is exact. 

Note that for $m = 1$ we thus obtain the split oscillator algebra. 
\end{example}

\begin{example} (a) If $\k = \spann \{x,y,c,d\}$ is the split oscillator algebra from 
Section~5.1, then $\k \cong T^*\g$ for the $2$-dimensional non-abelian 
subalgebra $\g := \spann\{x,d\}$ and 
Th.~3.2 in \cite{Bo97} implies that $(\k,\kappa_2)$ is isometrically isomorphic to $(T^*\g,\kappa)$, 
which provides another argument for the exactness of $\kappa_2$. 

(b) For $n = 3m$ and $\k = {\cal A}_{3m}$, we observe that the ideal 
$\fn := \spann\{T_{\frac{n+1}{2}}, \ldots, T_n\}$ is isotropic and abelian, so that 
\cite{Bo97}, Thm.~3.2 implies that $\k \cong T^*_\gamma(\g)$ for $\g := \k/\fn \cong 
{\cal A}_{\frac{n-1}{2}}$. Therefore Pelc's algebras provide other examples of 
exact twisted cotangent bundles. Using the canonical basis 
$(T_i)_{i = 0,\ldots, n}$ of 
$\k$ to define a section $\g \to \k$, we obtain 
$$ \gamma(T_i, T_j) = \left\{ 
  \begin{array}{cl} 
0 &\mbox{for $i + j \leq \frac{n-1}{2}$}  \\
\hat{i-j} \cdot T_{n-i-j}^* & \mbox{for $i + j > \frac{n-1}{2}$}
\end{array} \right. $$
In particular, $\gamma$ is non-zero. 
\end{example}

%%%%%%%%%%%%%%%%%%% new section %%%%%%%%%%%%%%%%%%%%%%%%%%%%%%
\section{The topological setting} \label{sec6}

In this section, we explain how the algebraic results from the preceding 
sections can be used in the topological setting. Actually these applications 
were our original motivation to study the work of Haddi and Zusmanovich. 

We now assume that $\K = \R$ or $\C$. 
Let $A$ be a unital commutative locally convex associative $\K$-algebra 
and $\k$ a locally convex $\K$-Lie algebra. 
We endow $\g = A \otimes \k$ 
with the projective tensor product topology, turning it into a locally convex 
space with the universal property that for each locally convex space $\z$ 
a bilinear map $\phi \: A \times \k \to \z$ is continuous if and only if the
corresponding linear map $\tilde \phi \: A \otimes \k \to \z$ is continuous.  
Then the Lie bracket on $\g$ is continuous because the quadrilinear map 
$$ A \times \fk \times A \times \fk \to A \otimes \fk, 
\quad (a,x,a',x') \mapsto aa' \otimes [x,x'] $$
is continuous and the continuous quadrilinear maps correspond to the 
continuous linear maps on \break $(A \otimes \k) \otimes (A \otimes \k)$. 

In the topological context, we consider for a locally convex space $\z$ 
the space $Z^2_c(\g,\z)$ of continuous cocycles and the 
subspace $B^2_c(\g,\z)$ of all coboundaries of the form $d_\g \ell$, where 
$\ell \: \g \to \z$ is a continuous linear map. In the topological context, 
the relation between the space 
$$ H^2_c(\g,\z) := Z^2_c(\g,\z)/B^2_c(\g,\z) $$
and the space of all linear maps from 
$$ H_{2,c}(\g) := Z_2(\g)/\oline{B_2(\g)} \to \z $$
is more complicated than in the algebraic setup (\cite{Ne02b}). 
To define the topological 
version of $H_2(\g)$, we have to use the closure of $B_2(\g)$ to obtain a 
Hausdorff topology on the quotient space. We always have a natural map 
$$ H^2_c(\g,\z) \to \Lin(H_{2,c}(\g),\z), $$
but in general there is no reason for this map to be injective or 
surjective. Therefore the homology space is much less interesting 
in the topological setting, and it often is easier to work directly with 
cocycles and coboundaries which is made possible by our results in Section~III. 

The flip involution on $\g \otimes \g$, endowed with the projective 
tensor product topology, is continuous, so that the kernel of the 
quotient map $\g \otimes \g \to \Lambda^2(\g), x \otimes y \mapsto x \wedge y$ 
is closed, which leads to a 
locally convex topology on $\Lambda^2(\g)$. Further the bracket map 
$b_\g \: \Lambda^2(\g) \to \g$ 
is continuous because it is induced from the continuous bracket map, 
which shows that its kernel $Z_2(\g)$ is closed. 

One easily verifies 
that the maps $p_\pm$ and $\sigma_\pm$ from Section~2 are continuous, and likewise 
that the maps 
$$ A \to A \otimes A, \quad a \mapsto a \otimes \1 \quad \mbox{ and } \quad 
 A \to S^2(A), \quad a \mapsto a \vee  \1 $$
are continuous. Therefore Lemma~\ref{lem2.1} yields a 
topological decomposition of the closed subspace $Z_2(\g)$ of $\Lambda^2(\g)$: 
$$ Z_2(\g) = (\Lambda^2(A) \otimes S^2(\k)) \oplus 
(A \otimes Z_2(\k)) \oplus (I_A \otimes \Lambda^2(\k)). $$
This implies that any continuous cocycle $f \: \g \times \g \to \z$ 
defines three continuous maps 
$$ f_1 \: \Lambda^2(A) \otimes S^2(\k) \to \z, \quad 
f_2 \: A \otimes \Lambda^2(\k) \to \z \quad \mbox{ and } \quad 
f_3 \: I_A \otimes \Lambda^2(\k) \to \z. $$
Conversely, three such continuous linear maps 
combine to a continuous $2$-cocycle of $\g$ if and only if 
they satisfy the conditions from Theorem~\ref{theo2}. 

If a continuous 
cocycle $f = f_1 +f_2 + f_3$ is contained 
in $B^2_c(\g,\z)$, then it vanishes on 
$Z_2(\g)$, which implies $f_1 = f_3 = 0$ and that 
$f_2$ is a continuous coboundary, i.e., there exists a continuous linear 
map $\ell \: \g \to \z$ with 
$$ f(ax, by) = f_2(ax,by) = \tilde f_2(ab)(x,y) 
= \ell(ab[x,y]) \quad \mbox{ for all } 
\quad a,b\in A, x,y \in \k. $$
Clearly, this implies that $\tilde f_2(A) \subeq B^2_c(\k,\z)$. 
If, conversely, $\tilde f_2(A) \subeq B^2_c(\k,\z)$, then there exists a 
linear map $h \: A \to \Lin(\k,\z)$ with 
$d_\k h(a) = \tilde f_2(a)$ for all $a \in A$, but it is not clear whether 
the corresponding
map $\tilde h \: A \times \k \to \z$ will be continuous. Therefore the 
exactness condition is quite subtle. 

If $\k$ is finite-dimensional, then the situation simplifies significantly. 
Then $B^2_c(\k,\z) = B^2(\k,\z)$ and if $\Lin(\cdot,\cdot)$ stands for 
``continuous linear maps'', then 
$$  \Lin(A \otimes \Lambda^2(\k),\z) 
\cong  \Lambda^2(\k)^* \otimes \Lin(A, \z), $$
 so that we may consider $f_2$ as a $2$-cocycle in $Z^2(\k, \Lin(A,\z))$.  
If this map vanishes on $B_2(\k)$, then there exists a linear map 
$h \: \k \to \Lin(A,\z)$ with 
$$ f_2(x,y)(a) = h([x,y])(a), \quad x,y \in \k, a \in A. $$
Then the map $\ell \: A \times \k \to \z, (a,x) \mapsto h(x)(a)$ 
is continuous and satisfies $f_2 = - d_\g \ell$. We thus get 
$$ B^2(\k,\Lin(A,\z)) \cong B^2_c(\g,\z). $$
We collect the previous remarks in the following theorem which is 
analoguous to Theorem~\ref{theo4}. It 
determines the structure of the second continuous cohomology space 
for current algebras. Let us denote by 
$\Omega_c^1(A)=J_A\,/\,\overline{J_A^2}$ 
the locally convex module of K\"ahler differentials for the locally convex 
commutative associative algebra $A$. 

\begin{theorem}
Let $\k$ be a finite dimensional Lie algebra over $\K$. Then the sequence 
$$ \{0\} \to H^2_c(\g/\g') \oplus \Lin(A, H^2(\k)) \ssmapright{\Phi} H^2_c(\g) 
\ssmapright{\Psi} \Lin((\Omega^1_c(A), \overline{d_A(A)}), (Z^3(\k)_\Gamma, B^3(\k)_\Gamma)) \to \{0\}$$ 
is exact.
\end{theorem}

\begin{proof} First we note that the short exact sequence 
$\0 \to \k' \to \k \to \k/\k' \to \0$ of finite-dimensional vector spaces splits. 
Since $\g' = A \otimes \k'$ is closed in $\g$, it follows that 
the short exact sequence 
$\0 \to \g' \to \g \to \g/\g' \to \0$ also splits topologically. 
As we have observed above, Theorem~\ref{theo3} and its corollaries remain true 
in the topological setting. For Corollary~\ref{cortheo2}, we use the topological 
splitting of $\k'$ in $\k$. We have also seen above that the 
corresponding description of the coboundaries 
remains valid because $\k$ is finite dimensional. 
Further, the topological splitting of $\g'$ implies that 
Lemma~\ref{lem4.1} remains true. This implies the injectivity of $\Phi$.

That $\ker\Psi=\im\Phi$ is shown as in the proof of Theorem~\ref{theo4}. 
Finally, the surjectivity of $\Psi$ follows from the fact that $f^{\flat}$ 
and $\beta$ can be chosen as continuous maps, because of 
the existence of (continuous) linear right inverses of surjective linear maps 
to finite dimensional vector spaces.
\end{proof}

\begin{example} We consider the special case where $M$ is a compact manifold and 
$A = C^\infty(M,\R)$ the Fr\'echet algebra of all smooth real-valued functions on~$M$.
According to \cite{Mai} or \cite{Co85}, the universal topological differential module 
of $A$ is given by $\Omega^1_c(A) \cong \Omega^1(M,\R)$, the space of smooth $\R$-valued 
$1$-forms on $M$, and the de Rham-differential $d \: C^\infty(M,\R) \to \Omega^1(M,\R)$ 
is a universal continuous derivation. It follows in particular that 
the space $d_A(A)$ is the space of exact $1$-forms, which is non-zero. 

Now let $\k$ be a finite-dimensional real Lie algebra and 
$$\g := A \otimes \k \cong C^\infty(M,\k). $$
Up to cocycles vanishing on $\g \times \g'$, all continuous cohomology classes 
in $H^2_c(\g)$ are then represented by sums $f = f_1 + f_2$, where 
$$ \tilde f_1 \: A \times A \to \Sym^2(\k)^\k $$
is an alternating continuous linear map for which there is a continuous 
linear map 
$$ f_1^\flat \: \Omega^1(M,\R) \to \Sym^2(\k)^\k 
\quad \mbox{ with } \quad 
\tilde f_1(a,b) = f_1^\flat(a \cdot d(b) - b \cdot d(a)), $$
and 
$$ d_\k(\tilde f_2(a)) = -\Gamma(f_1^\flat(da)) \quad \mbox{ for all} \quad a \in A. $$

We interprete the continuous linear map $f_1^\flat$ as a $\Sym^2(\k)^\k$-valued 
{\it current on $M$}. It is a closed current if and only if it vanishes on exact forms. 
Typical examples of such currents arise 
from pairs $(\xi, \kappa)$, where $\xi \: [0,1] \to M$ is a piecewise smooth path and 
$\kappa \in \Sym^2(\k)^\k_{\rm ex}$ via 
$$ f_1^\flat(\alpha) := \Big( \int_\xi \alpha\Big) \cdot \kappa, $$
but these examples satisfy $\gamma \circ f_1^\flat = 0$. 
%{\bf Is there a geometric way to construct more non-trivial examples???} 
\end{example}

%%%%%%%%%%%%%%%%%%% new section %%%%%%%%%%%%%%%%%%%%%%%%%%%%%%
\section{Appendix: A useful exact sequence} \label{sec7}

The following section is very much based on information and hints we got from 
M.~Bordemann (\cite{Bo97}). 

\begin{definition}
Let $\k$ be a Lie algebra and $\fa$ a $\k$-module. 
We denote the action as $\k \times \fa \to \fa$ by $(x,a) \mapsto x.a$. 
On the space $C^p(\k,\fa)$ of $\fa$-valued Lie algebra cochains we have  a natural action of $\k$ 
denoted by 
$$(x.\omega)(x_1,\ldots, x_p) 
= x.\omega(x_1,\ldots, x_p) 
- \sum_{i=1}^p \omega(x_1, \ldots, x_{i-1}, [x,x_i], x_{i+1}, \ldots, x_p). $$

For $p,q \in \N_0$, we consider the injection  
$$ \tilde T_p \: C^{p+q}(\k,\fa) \to C^p(\k,C^q(\k,\fa)), \quad 
(\tilde T_p f)(x_1,\ldots, x_p)(y_1, \ldots, y_q) := f(x_1, \ldots, x_p, y_1, \ldots, y_q). $$
{}From the action of $\k$ on the 
spaces $C^q(\k,\fa)$, we obtain 
Lie algebra differentials 
$$ d_\k' \: C^p(\k,C^q(\k,\fa)) \to C^{p+1}(\k,C^q(\k,\fa)) $$
and we also have 
$$ d_\k'' \: C^p(\k,C^q(\k,\fa)) \to C^p(\k,C^{q+1}(\k,\fa)), \quad 
\omega \mapsto d_\k \circ \omega $$
satisfying on $C^{p+q}(\k,\fa)$  the identity 
\begin{eqnarray}
  \label{eq:homform}
\tilde T_{p+1} \circ d_\k = d_\k' \circ \tilde T_{p} + (-1)^{p+1}  d_\k'' \circ \tilde T_{p+1}. 
  \end{eqnarray}
(cf.\ \cite{HS53}, Lemma~1). 
\end{definition} 

Specializing to the trivial module $\fa = \K$, we obtain in particular the maps 
$$ \tilde\alpha_p := \tilde T_{p-1} \: C^p(\k) \to C^{p-1}(\k,C^1(\k)) = C^{p-1}(\k,\k^*), $$
which, in view of equation (\ref{eq:homform}), 
commute with the respective Lie algebra differentials 
because $d_\k'' \circ \tilde T_p$ vanishes on $C^p(\k,\K)$. Hence they induce linear maps 
$$ \alpha_p \: H^p(\k) \to H^{p-1}(\k,\k^*), \quad [\omega] \mapsto [\tilde \alpha_p(\omega)]. $$

For the $\k$-module $\k^*$, the subspace $d_\k \k^*$ of 
$C^1(\k,\k^*)$ consists of maps 
whose associated bilinear map is alternating. We thus have a well-defined 
map 
$$ S \: C^1(\k,\k^*)/B^1(\k,\k^*) \to \Sym^2(\k) $$
which is a morphism of $\k$-modules. We now obtain maps 
$$ \tilde\beta_p = S \circ \tilde T_{p-1}\: C^p(\k,\k^*) \to C^{p-1}(\k, \Sym^2(\k)) $$
satisfying 
\begin{eqnarray}
\tilde\beta_p \circ d_\k 
&=& S \circ \tilde T_{p-1} \circ d_\k 
= S \circ (d_\k' \circ \tilde T_{p-2}) 
= d_\k' \circ S \circ \tilde T_{p-2} = d_\k' \circ \tilde \beta_{p-1}.
\end{eqnarray}
Hence $\tilde\beta_{p}$ induces a linear map 
$$ \beta_{p} \: H^p(\k,\k^*) \to H^{p-1}(\k,\Sym^2(\k)). $$
{}From the construction, we immediately get $\tilde\beta_p \circ \tilde\alpha_{p+1} = 0$, which 
leads to $\beta_p  \circ \alpha_{p+1} =0$. 

\begin{proposition} \label{transfer-seq} 
For any Lie algebra $\k$, we obtain with $\gamma(\kappa) := [\Gamma(\kappa)]$ 
an exact sequence 
$$ \{0\} \to H^2(\k) \ssmapright{\alpha_2} 
H^1(\k,\k^*) \ssmapright{\beta_1} \Sym^2(\k)^\k 
\ssmapright{\gamma} H^3(\k) 
\ssmapright{\alpha_3} H^2(\k,\k^*) \ssmapright{\beta_2} H^1(\k,\Sym^2(\k)).  $$
\end{proposition}

\begin{proof}
To see that for each cocycle $\omega \in Z^1(\k,\k^*)$ the symmetric bilinear form 
$\tilde\beta_1(\omega)$ is invariant, we note that 
$$ \tilde\beta_1(\omega)([x,y],z) = \omega([x,y])(z) +\omega(z)([x,y]), $$
and if $\omega$ is a cocycle, this can be written as 
$$ \tilde\beta_1(\omega)([x,y],z) = (x.\omega(y))(z) - (y.\omega(x))(z) + \omega(z)([x,y]) 
= \omega(y)([z,x]) + \omega(x)([y,z]) + \omega(z)([x,y]), $$
showing that this trilinear form is alternating, and hence that $\tilde\beta_1(\omega)$ is 
invariant. 

{\bf Exactness in $H^2(\k)$:} We only have to show that $\alpha_2$ is injective. 
If $\omega \in Z^2(\k)$ satisfies $\tilde\alpha_2(\omega) = d_\k \eta$ for some 
$\eta \in \k^*$, then 
$$ \omega(x,y) = (d_\k\eta)(x)(y) = (x.\eta)(y) = - \eta([x,y]), $$
which implies that $\omega$ is a $2$-coboundary. 

{\bf Exactness in $H^1(\k,\k^*)$:} Clearly $\beta_1 \circ \alpha_2=0$. 
If, conversely, $\beta_1([\omega]) = 0$, then $\omega \: \k \to \k^*$ is a linear map 
whose associated bilinear form $\tilde\omega(x,y) := \omega(x)(y)$ 
is alternating. In this situation, we have 
\begin{eqnarray}
d\tilde\omega(x,y,z) 
&=& -\omega([x,y])(z) - \omega([y,z])(x) - \omega([z,x])(y)
= -\omega([x,y])(z) + \omega(x)([y,z]) + \omega(y)([z,x]) \cr
&=& \big(- \omega([x,y]) - y.\omega(x) + x.\omega(y)\big)(z)
=  (d_\k \omega)(x,y)(z).
\end{eqnarray} 
We conclude that $\tilde\omega$ is a cocycle if and only if $\omega$ is one, and from that 
we derive that $\ker \beta_1 = \im \alpha_2$. 

{\bf Exactness in $\Sym^2(\k)^\k$:} 
Next we show that $\gamma \circ \beta_1 = 0$. So let $\omega \in Z^1(\k,\k^*)$ 
and write $\tilde\omega = \omega_+ + \omega_-$, where $\omega_+$ is 
symmetric and $\omega_-$ is alternating. 
Then 
\begin{eqnarray}
\Gamma(\tilde\beta_1(\omega))(x,y,z) 
&=& \omega([x,y])(z) + \omega(z)([x,y])
= \omega(y)([z,x]) + \omega(x)([y,z]) + \omega(z)([x,y]), 
\end{eqnarray}
and the closedness of $\omega$ also shows that 
$\sum_\cyc \omega([x,y])(z) = 2 \sum_\cyc \omega(x)([y,z]),$
which leads to 
$$ \Gamma(\tilde\beta_1(\omega))(x,y,z) 
= \sum_\cyc \omega([x,y])(z) - \omega(z)([x,y]) 
= 2 \sum_\cyc \omega_-([x,y],z) = - 2 d_\k \omega_-(x,y,z). $$
Hence $\Gamma(\tilde\beta_1(\omega))$ is always exact, so that $\gamma \circ 
\beta_1$ vanishes on the level of cohomology spaces. 

To see that $\ker \gamma \subeq \im \beta_1$, suppose that $\kappa$ is an exact invariant 
bilinear form and $\eta \in C^2(\k)$ satisfies $d\eta = - \Gamma(\kappa)$. 
Then 
$\omega(x)(y) := \kappa(x,y) + \eta(x,y)$
defines a linear map $\omega \: \k \to \k^*$ with 
\begin{eqnarray}
(d_\k\omega)(x,y)(z) 
&=& (x.\omega(y) - y.\omega(x) - \omega([x,y]))(z) \cr
&=& \omega(y)([z,x]) + \omega(x)([y,z]) - \omega([x,y])(z) \cr
&=& \kappa([x,y],z) + \eta(y,[z,x]) + \eta(x,[y,z]) - \eta([x,y],z) \cr
&=& \kappa([x,y],z) + d_\k\eta(x,y,z) = (\Gamma(\kappa) + d_\k\eta)(x,y,z) = 0.
\end{eqnarray}

{}From the preceding calculation, we also see by putting $\eta = 0$, that 
the linear map $\tilde\kappa \: \k \to \k^*$ defined by an invariant symmetric bilinear form $\kappa$ 
is a $1$-cocycle if and only if $\Gamma(\kappa)$ vanishes. 

{\bf Exactness in $H^3(\k)$:} 
The transfer formula for differentials implies 
that an alternating trilinear form $\omega$ on $\k$ is a $3$-cocycle if and only if 
the corresponding alternating bilinear form $\tilde\alpha_3(\omega)$ is a $2$-cocycle. 
Therefore the image of $\alpha_3$ consists of those cohomology classes having a representing 
cocycle whose associated trilinear form is alternating. 

For $\kappa \in \Sym^2(\k)^\k$, the corresponding $3$-cocycle 
$\Gamma(\kappa)$ and the corresponding linear map 
$\tilde\kappa \: \k \to \k^*$, we have  
$$\tilde\alpha_3(\Gamma(\kappa))(x,y) = \kappa([x,y],\cdot) = - (d\tilde\kappa)(x)(y) $$
because 
\begin{eqnarray}
d\tilde\kappa(x,y)(z) 
&=& (x.\tilde\kappa(y))(z) - (y.\tilde\kappa(x))(z) - \tilde\kappa([x,y])(z) \cr
&=& -\kappa(y, [x,z]) + \kappa(x,[y,z]) - \kappa([x,y],z) 
= \kappa([x,y],z). 
\end{eqnarray}
We conclude that $\tilde\alpha_3(\Gamma(\kappa))$ is exact, so that 
$\tilde\alpha_3 \circ \Gamma$ induces the trivial map $\alpha_3 \circ \gamma \: 
\Sym^2(\k)^\k \to H^2(\k,\k^*)$. 

Let $f \in C^1(\k,\k^*)$ and write $\tilde f(a,b) = f(a)(b)$. We then have 
\begin{eqnarray}
df(x,y)(z) 
&=& (x.f(y) - y.f(x) - f([x,y]))(z)
= f(y)([z,x]) + f(x)([y,z]) - f([x,y])(z) \cr
&=& f(y)([z,x]) - (y.\tilde f)(x,z). 
\end{eqnarray}
This map is alternating in $(x,y)$, and it is alternating in $(x,z)$ if and only if 
$y.\tilde f$ is alternating. Writing $\tilde f = \tilde f_+ + \tilde f_-$ for the decomposition 
of $\tilde f$ into symmetric and alternating components, this is equivalent to 
$y.\tilde f_+ = 0$. We conclude that $df(x,y)(z)$ is alternating if and only if 
$\tilde f_+$ is invariant. 

To verify the exactness in $H^3(\k)$, we now 
assume that $\omega \in Z^3(\k)$ satisfies $\alpha_3(\omega) \in B^2(\k,\k^*)$, 
i.e., $\alpha_3([\omega]) = 0$. Then there exists an $f \in C^1(\k,\k^*)$ with 
$\tilde\alpha_3(\omega) = d_\k f$, and the preceding paragraph implies that 
$\tilde f_+$ is an invariant symmetric bilinear form on $\k$ satisfying 
$$ \tilde\alpha_3(\omega) =   d_\k f = - \Gamma(\tilde f_+) + d_\k f_-, $$
where $f = f_+ + f_-$ corresponds to the decomposition $\tilde f= \tilde f_+ + \tilde f_-$. 
We conclude that $[\tilde\alpha_3(\omega)] = -[\Gamma(\tilde f_+)]$, which implies 
exactness in $H^3(\k)$. 

{\bf Exactness in $H^2(\k,\k^*)$:} 
We claim that $\ker \beta_2 = \im \alpha_3$. To verify this claim, 
pick $\omega \in Z^2(\k,\k^*)$ for which $\tilde\beta_2(\omega)$ is exact, i.e., 
there exists a symmetric bilinear form $\kappa \in \Sym^2(\k)$ with 
$\tilde\beta_2(\omega) = d_\k \kappa$, i.e., 
for $x,y,z\in\k$ we have 
$$ \omega(x,y)(z) + \omega(x,z)(y) = (x.\kappa)(y,z) 
= - \kappa([x,y],z) - \kappa(y,[x,z]). $$
Let $\tilde\eta\in C^1(\k,\k^*)$ and write $\eta$ for the corresponding bilinear map on $\k$ 
with $\eta(x,y) = \tilde\eta(x)(y)$. 
Then 
$$ d \tilde\eta(x,y)(z) = (x.\tilde\eta(y) - y.\tilde\eta(x) - \tilde\eta([x,y]))(z) 
= -\eta(y,[x,z]) + \eta(x,[y,z]) - \eta([x,y],z). $$
Therefore 
\begin{eqnarray}
\tilde\beta_2(d_\k\tilde\eta)(x)(y,z)  
&=& \tilde\eta(y)([z,x]) - \tilde\eta([x,y])(z) + \tilde\eta(z)([y,x]) - \tilde\eta([x,z],y) \cr
&=& 2 (\tilde\eta_+([y,x])(z) + \tilde\eta_+([z,x])(y)) 
= 2 (x.\eta_+)(y,z),
\end{eqnarray}
and this leads to 
$$ \tilde\beta_2(d_\k\tilde\eta) = 2 d_\k \eta_+ = d_\k(\tilde\beta_1(\eta)). $$
Since $\tilde\beta_1(C^1(\k,\k^*)) = \Sym^2(\k)$, we find some $\tilde\eta \in C^1(\k,\k^*)$ with 
$\tilde\beta_1(\tilde\eta) = \kappa$, and then 
$$ \tilde\beta_2(\omega - d_\k \eta) = d_\k \kappa - d_\k \eta_+ = 0, $$
so that for $\omega' := \omega - d_\k \eta \in Z^2(\k,\k^*)$ the corresponding trilinear map 
is alternating. This means that $[\omega] = [\omega'] \in \im \alpha_3$. 
\end{proof}

\end{document}